\documentclass[a4paper]{article}
\usepackage[english]{babel}
\usepackage[utf8]{inputenc}
\usepackage{amsmath,amsthm}
\usepackage{amsfonts,amssymb}

\usepackage{mathtools}

\usepackage{url}

\usepackage{tikz,tikz-cd} 
 \usetikzlibrary{decorations.pathmorphing} 
\usetikzlibrary{babel}

\usepackage{xifthen} 

\newcommand{\N}{\mathbb{N}}

\newcommand{\G}{\mathbb{G}}

\newcommand{\Pb}{\mathbb{P}}
\newcommand{\Ob}{\mathbb{O}}

\newcommand{\Ecal}{\mathcal{E}}

\newcommand{\Ocal}{\mathcal{O}}

\newcommand{\Dcal}{\mathcal{D}} 
\newcommand{\Fcal}{\mathcal{F}}

\newcommand{\Ccal}{\mathcal{C}}

\newcommand{\Glob}{\mathbb{G}}

\usepackage{titlesec}
\titleformat{\subsubsection}[runin]{\normalfont}{\thesubsubsection}{0pt}{}[.]

\renewcommand{\thesubsubsection}{\arabic{section}.\arabic{subsection}.\arabic{subsubsection}}
\csname @addtoreset\endcsname{subsubsection}{section}

\newcommand{\block}[1]
{

\par \subsubsection{} #1

\bigskip}

\newcommand{\blockn}[1]{\par #1 \bigskip}

\newcommand{\Th}[1]
	{
	\bigskip	
	\textbf{Theorem : }{\itshape #1}
		
	\bigskip
	}

\newcommand{\Prop}[1]
	{

	\bigskip
	
	\textbf{Proposition : }{\itshape #1}
		
	\bigskip
	
	}

\newcommand{\Conjecture}[1]
	{

	\bigskip
	
	\textbf{Conjecture : }{\itshape #1}
		
	\bigskip
	
	}

\newcommand{\Cor}[1]
	{

	\bigskip
	
	\textbf{Corollary : }{\itshape #1}	
		
	\bigskip

	}

\newcommand{\Lem}[1]
	{

	\bigskip
	
	\textbf{Lemma : }{\itshape #1}
		
	\bigskip
	
	}

\newcommand{\Def}[1]
	{
	
	\bigskip
	
	\textbf{Definition : }{\itshape #1}
	
	\bigskip
	
	}

\newcommand{\Dem}[1]{
	
	\smallskip
	
	\textbf{Proof : } \par
	 {#1} $\square$
	 
	 \bigskip
}

\begin{document}

\pagestyle{plain}
\title{Non-unital polygraphs form a presheaf category}

\author{Simon Henry}
\date{}

\maketitle

\begin{abstract}

We prove, as claimed by A.Carboni and P.T.Johnstone, that the category of non-unital polygraphs, i.e. polygraphs where the source and target of each generator are not identity arrows, is a presheaf category. More generally we develop a new criterion for proving that certain classes of polygraphs are presheaf categories. This criterion also applies to the larger class of polygraphs where only the source of each generator is not an identity, and to the class of ``many-to-one polygraphs'', producing a new, more direct, proof that this is a presheaf category. The criterion itself seems to be extendable to more general type of operads over possibly different combinatorics, but we leave this question for future work.

In an appendix we explain why this result is relevant if one wants to fix the arguments of a famous paper of M.Kapranov and V.Voevodsky and make them into a proof of C.Simpson's semi-strictification conjecture. We present a program aiming at proving this conjecture, which will be continued in subsequent papers.
\end{abstract}

\tableofcontents

\renewcommand{\thefootnote}{\fnsymbol{footnote}} 
\footnotetext{This work was supported by the Grant agency of the Czech republic under the grant P201/12/G028.}
\renewcommand{\thefootnote}{\arabic{footnote}} 

\section*{Introduction}

\blockn{Some time ago I started studying the famous false proof by M.Kapranov and V.Voevodsky in \cite{kapranov1991infty} that homotopy types can be represented by strict $\infty$-categories whose arrows are weakly invertible, hopping to prove the conjecture by C.Simpson in \cite{simpson1998homotopy} that this false proof can be made into a correct proof of the fact that every homotopy types can be represented as an $\infty$-groupoid whose associativity and exchange rule hold strictly, but where units and inverses are weak.

This analysis lead me to the conclusion that the main problem with \cite{kapranov1991infty} is in the choice of the category of higher categorical diagrams they are using. A more detailed analysis showed that constructing a category of diagrams having the correct properties for this proof to work appears to be closely related to another problem: Showing that certain classes of polygraphs are presheaf categories. The present paper is devoted to this second problem, but we have included in an appendix \ref{App-Simpson-voevodsky-Kapranov} a presentation of our ideas on why we think the argument of M.Kapranov and V.Voevodsky fails and how the results of the present paper can rescue their ideas and, maybe, lead to a proof of C.Simpson conjecture. In this appendix we also explain precisely why we have not succeeded in proving this conjecture yet, and what remains to be done. This line of work will be pursued further in \cite{henry2018regular}, leading to a proof of a form of C.Simpson conjecture. To some extent, this appendix constitute an introduction to the present paper in the sense that it explains its main motivations, but the paper can be read completely independently of this appendix and the results of the paper are interesting on their own right for the theory of polygraphs. Also, some understanding of the theory developed in the present paper will be needed to follow the discussion at the end of the appendix, which is why we postponed it to the end of the paper.}

\blockn{Our main result is that the category of ``non-unital'' polygraphs, or ``positive polygraphs'' as we will call them, i.e. polygraphs where the source and targets of each generator is not an identity arrow is a presheaf category.

This was claimed without a proof by P.T.Johnstone and A.Carboni in \cite{carboni2004corrigenda} after they noted that, because of the Eckmann-Hilton argument, their claim in \cite{carboni1995connected} that the category of all polygraphs is a presheaf category was false. We prove more generally that the category of polygraphs where the source (or equivalently the target) of each generator is not an identity arrow is a presheaf category, and we prove that subclasses of polygraphs (in the sense of definition \ref{Def_classofPG}) of this class are all presheaf categories.

\bigskip

This in particular gives the first direct proof that the category of ``many-to-one'' polygraphs, i.e. polygraphs where the target of each generator is a generator, is a presheaf category. This is known due to an indirect proof that this category is equivalent to the category of opetopic sets, following from \cite{hermida2000weak} and \cite{harnik2008computads} together. Since the online publication of the first version of this paper, a direct proof of this equivalence has also been given in \cite{thanh2018equivalence}.

\bigskip

We also believe that the methods used to obtain this result are of independent interest and might be applied to other types of polygraphs, for examples polygraphs corresponding to less strict notion of higher categories, or to non-globular structure. We hope to come back to this in a future work.

\bigskip

We postpone the introduction of the main ideas involved in that proof to the beginning of section \ref{section_goodclass}, after a short introduction to strict $\infty$-categories and polygraphs in section \ref{section_preliminaries}. 
}

\section{Polygraphic preliminaries}
\label{section_preliminaries}
\subsection{$\infty$-categories}

\block{\label{Def_StreetInfinityCat}We will use R.Street ``one type'' definition of $\infty$-category (from \cite{street1987algebra}, see also \cite{steiner2004omega}): 

\Def{An $\infty$-category is a set $X$ together with unary operations $\pi^-_k$ and $\pi^+_k$ called respectively ``$k$-source'' and ``$k$-target'' for $k\geqslant 0$ and partially defined binary operations $\#_n$ called ``$n$-composition'' for $n\geqslant 0$ satisfying the following axioms:

\begin{enumerate}

\item $ x \#_n y $ is defined if and only if $\pi^+_n(x)=\pi^-_n(y)$.

\item For every $x \in X$ there exists an $n$ such that $\pi^-_n(x)=\pi^+_n(x)=x$.

\item For any $x \in X$, one has:

\[ \pi^{\epsilon}_n \pi^{\delta}_m x = \left\lbrace \begin{array}{l l} \pi^{\epsilon}_n x & \text{if } n <m \\ \pi^{\delta}_m x & \text{if } n \geqslant m\end{array} \right. \]

where $\epsilon$ and $\delta$ are arbitrary signs.

\item For any $x \in X$ one has $\pi^-_n(x) \#_n x = x \#_n \pi^+_n(x) = x $.

\item For all $x,y$ such that the $x \#_n y$ is defined,

\[ \pi^-_n(x \#_n y) = \pi^-_n(x) \]
\[ \pi^+_n(x \#_n y) = \pi^+_n(y) \]

And if $k>n$:

\[ \pi^{\epsilon}_k(x \#_n y) = \pi^{\epsilon}_k(x) \#_n \pi^{\epsilon}_k(y)\]

\item $x \#_n (y \#_n z) = (x \#_n y) \#_n z$ when either side is defined.

\item If $k <n$

\[ (x \#_n y) \#_k ( z \#_n w) = (x \#_k z) \#_n (y \#_k w) \]

when the left hand side is defined.

\end{enumerate}

A morphism of $\infty$-categories is a function commuting to all the structural functions ($\pi^{\epsilon}_k$ and $\#_k$ ).

}

A few remarks in order to clarify the intended meaning of the definition and to relate it to the ``globular'' definition (i.e. as globular set with operations):

\begin{itemize}

\item In terms of the ``globular'' definition of $\infty$-categories, the underlying set $X$ corresponds to the increasing union of the sets of $n$-arrows for all $n$, where an $n$-arrows is seen as a $(n+1)$-arrow by identifying it with its identity endomorphism.

\item $\pi^-_n(x)$ and $\pi^+_n(x)$ denotes respectively the ``$n$-dimensional source and target'' of $x$, i.e. the first iterated source/target of $x$ which is an $n$-arrow. In particular the set of $n$-arrows is exactly the set of $x \in X$ such that $\pi^{-}_n x=x$ (or equivalently $\pi_n^+ x = x$) and if $x$ is an $n$-arrow then its usual source and target are the $(n-1)$-arrow $\pi^-_{n-1}(x)$ and $\pi^+_{n-1}(x)$.

\item The operation $x \#_n y$ corresponds to the composition of two arrows along a common boundary of dimension $n$, note that for $\infty$-categories we use the ``reverse'' or ``diagrammatic'' composition order: $x \#_n y$ is defined when the $n$-target of $x$ is the $n$-source of $y$. The difference of convention for composition order between $\infty$-categories and all other categories will not be a problem as we will always consider strict $\infty$-categories as just algebraic structures of interest, and not as actual categories (for example ``equivalences'' of $\infty$-categories will play no role in the present paper, only isomorphisms).

\item Note that this operation $x \#_n y$ allows to compose arrows of different dimensions. In terms of the globular definition, if $x$ is of dimension lower than $y$ then $x \#_n y$ corresponds to the $n$-composition of $y$ with the iterated identity of $x$ of same dimension as $y$. This will become very important when we discuss non-unital $\infty$-categories in \ref{Def_nonUnitalInfinityCat}.

\item We have not given rules for the value of $\pi^{\epsilon}_k(x \#_n y)$ when $k<n$ but it follows from the axiom for $\pi^{\epsilon}_n(x \#_n y)$ that $\pi^{\epsilon}_k(x \#_n y)= \pi^{\epsilon}_k(x) = \pi^{\epsilon}_k(y)$ when $k<n$.

\end{itemize}

}

\block{If $X$ is an $\infty$-category, elements of $X$ will be called the arrows of $X$. An $n$-arrow of $X$ is an arrow $x \in X$ such that $\pi^{\epsilon}_n(x)=x$ (it holds for $\epsilon=+$ if and only if it holds for $\epsilon=-$). For any arrow $x$, the $n$-source and $n$-target $\pi^{-}_n x $ and $\pi^+_n(x)$ are $n$-arrows. Any $n$-arrow can also be seen as an $m$-arrow for $m>n$ and these corresponds to identity arrows.

\bigskip

The dimension of an arrow $x$ is the smallest $n$ such that $x$ is an $n$-arrow. In particular an arrow of dimension $n$ is not the same as an $n$-arrow: an arrow of dimension $n$ is exactly a non-identity $n$-arrow. 

\bigskip

Two $n$-arrows $x$ and $y$ are said to be parallel if $\forall \epsilon, \pi^{\epsilon}_{n-1}x = \pi^{\epsilon}_{n-1} y$. An $(n+1)$-arrow $z$ is said to be ``from $x$ to $y$'', written $z: x \rightarrow y$ if $\pi^-_n z = x$ and $\pi^+_n z =y$, this implies that $x$ and $y$ are parallel.
}

\subsection{Polygraphs}
\label{subsection_polygraphs}

\blockn{Polygraphs were first introduced under the name ``computads'' by R.Street in \cite{street1976limits} in the framework of $2$-categories. The general $n$-categorical notion of computad is hinted at, also by R.Street, in \cite{street1987algebra} and as far as I know first appears explicitly in \cite{power1991n} and in \cite{burroni1993higher}. The name of polygraphs is due to A.Burroni, when he re-introduced the notion independently in \cite{burroni1993higher}.

\bigskip

In this subsection we will just review the definition and basic property of polygraphs, without any proof.

}

\blockn{The definition of the notion of polygraph is done by induction on the dimension, together with the definition of the ``free $\infty$-category'' $X^*$ on a polygraph $X$.}

\block{\Def{\begin{itemize}

\item A $0$-polygraph is just a set. The free $\infty$-category on a $0$-polygraph $X$ is the discrete $\infty$-category with only $0$-dimensional arrows given by the elements of $X$. The elements of the set $X$ are called the $0$-cells of $X$.

\item An $(n+1)$-polygraph $X$ is the data of an $n$-polygraph $X_n$ together with a set of so-called $(n+1)$-cells of $X$ and two maps $\pi^-_n$ and $\pi^+_n$ which attach to every $(n+1)$-cell a pair of parallel $n$-arrows in $(X_n)^*$. The free $\infty$-category generated by an $(n+1)$-polygraph is the $\infty$-category obtained by freely adding an $(n+1)$-arrow to $(X_n)^{*}$ for each $(n+1)$-cell of $X$ (with the source and target given by the $\pi^{\epsilon}_n$), i.e. it is defined by the following universal property:

A morphism $f$ from $X^*$ to another $\infty$-category $Z$ is the data of a morphism $f_n :(X_n)^* \rightarrow Z$ together with the choice of an $(n+1)$-arrow $f(x)$ in $Z$ for all $(n+1)$-cell $x$ of $X$ such that $\pi^{\epsilon}_n(f(x))=f_n(\pi^{\epsilon}_n(x))$.

\item A polygraph $X$ is a tower of, for each $n \in \N$, an $n$-polygraph $X_n$ where for all $n$ the underlying $n$-polygraph of $X_{n+1}$ is $X_n$. The free infinity category $X^*$ generated by a polygraph $X$ is defined as the inductive limit of the $\infty$-categories $X_n^*$.

\end{itemize}

}

For example, a $1$-polygraph is just an ordinary oriented graph, and the free $\infty$-category generated by it is the usual free category generated by a graph, i.e. the category of paths in the graph. In a $2$-polygraph one has additionally some $2$-cells, each $2$-cell being a ``$2$-arrow'' between two (oriented) path in the graph with the same source and target.

An $n$-polygraph will always be considered as an $(n+1)$-polygraph with no $(n+1)$-cells, and more generally, as a polygraph with no $k$-cells for $k>n$. This is compatible with the construction of the free $\infty$-category.

The notion of morphisms of polygraphs is also defined by induction, and such that the free $\infty$-category construction becomes a functor: A morphism of $0$-polygraphs is just a function between sets, which can also be seen as a functor between the corresponding $\infty$-categories. A morphism $f$ between two $(n+1)$-polygraphs is the data of a morphism $f_n$ between their underlying $n$-polygraphs together with a function $f$ between their sets of $(n+1)$-cells which are compatible to the maps $\pi^{\epsilon}_n$ in the sense that $\pi^{\epsilon}_n(f(x))= (f_n)^*(\pi^{\epsilon}_n(x))$. Such a data clearly induces a morphism between the corresponding free $\infty$-categories.

This gives a category of $n$-polygraphs, and more generally a category of polygraphs. The free $\infty$-category construction is a functor from the category of polygraphs to the category of $\infty$-categories.

The category of all polygraphs is denoted by $\Pb$, the (full) subcategory of $n$-polygraphs by $\Pb_n$.

}

\block{If $X$ is a polygraph, each $n$-cell $x$ of $X$ defines an $n$-arrow of $X^*$, also denoted $x$. It is in fact always an $n$-dimensional arrow of $X^*$ as the $\pi^{\epsilon}_{n-1}(x)$ are elements of $X_{n-1}$ and so cannot be equal to $x$.

These specific arrows of $X^*$ will be called the generators of $X^*$. A generator of $X^*$ of dimension $n$ is exactly the same as a $n$-cell of $X$. We want to avoid as much as we can to use the word ``cell'' for anything related to $X^*$ to avoid the confusion between cells of $X$ and the possible use of this word for the arrows of an $n$-category.

\bigskip

}

\block{\label{Discus:SyntacticalPropOfPolygraphs}In general, the set of all arrows of $X^*$ admit a description as the set of expressions formed using the composition operations $(\_ \#_k \_)$ and the cells of $X$, which are well formed (syntactically) and well defined (i.e. such that the composition involved are well defined) up to an equivalence relation generated by the axioms of $\infty$-categories. A precise formulation of this statement can be found in \cite{penon1999approche} or in \cite{makkai2005word}. We will not really use this explicitly.

What will be more important for us are the following simpler properties, which are all easy consequences of the above claim, but can all be proved more directly:

\begin{itemize}

\item The generators of $X^*$, i.e. the cells of $X$, are exactly the arrows that cannot be written in a non-trivial way as a composite in $X$. In particular, if an $\infty$-category is free over a polygraph, then the polygraph is uniquely determined up to unique isomorphism. Morphisms of polygraphs from $X$ to $Y$ are the same as morphisms of $\infty$-categories $f:X^* \rightarrow Y^*$ which send any generator of $X^*$ to a generator of $Y^*$. In particular, any isomorphism between $X^*$ and $Y^*$ comes from an isomorphism of polygraphs between $X$ and $Y$.

\item One can prove a property by ``induction on the arrows of a polygraph'': Let $P$ be a property of arrows of a polygraph $X$. In order to prove that $P(x)$ holds for all $x$, it is enough to prove the following: $P(a)$ holds for all $0$-dimensional generators $a$; If $P(x)$ holds for all $n$-arrows of $X^*$ for some $n$, then $P(a)$ holds for all $(n+1)$-dimensional generators; If $P(x)$ and $P(y)$ both holds for $x$ and $y$ two $k$-composable $n$-arrows and $P(z)$ holds for all $(n-1)$-arrows then $P(x \#_k y)$ holds.

\item If $a$ is an $n$-arrow of $X^*$ and $b$ is an $n$-cell of $X$, the ``number of times $b$ appears in any expression of $a$'' is a well defined number which does not depend on the expression\footnote{This only applies to expressions that do not contain $\#_m$ for $m \geqslant n$. Indeed these compositions can always be eliminated using the unit law (point $4.$ of \ref{Def_StreetInfinityCat}) but doing so can change the number of occurrences of some generators.} of $a$ in terms of the generators. Moreover $a$ is an identity $n$-arrow (i.e. an $(n-1)$-arrow) if and only if none of the generators of dimension $n$ appears in $a$. This can be deduced from the remark above and the fact that none of the axioms of $\infty$-category change that number of appearances, or more directly using the ``linearization'' functor (sending an $\infty$-category to a chain complex it generates) which implement this counting function (see subsection $4.3$ of \cite{metayer2008cofibrant}). Composition is additive with respect to this counting function, i.e. the number of times $a$ appears in $x \#_k y$ is the sum of the number of time it appears in $x$ and in $y$. It does not seem possible to define such a number when the dimension of $a$ and $x$ are different (a counting function of this kind is defined in \cite{makkai2005word}, but it does not really corresponds to a number of appearances).

\end{itemize}

So to some extent, polygraphs are just particular $\infty$-categories, which are ``free''. In fact it has been proved by F.Metayer in \cite{metayer2008cofibrant} that polygraphs are exactly the cofibrant objects of the folk model structure on $\infty$-categories.

If $X$ and $Y$ are polygraphs, a morphism from $X^*$ to $Y^*$ will be said to be polygraphic if it send generators to generators. As mentioned above, polygraphic morphisms are exactly the same thing as morphisms of polygraphs, and are in general more restrictive than morphisms of $\infty$-categories between $X^*$ and $Y^*$.
}

\block{\label{RightAdjointPG}The functor $X \rightarrow X^*$ which sends a polygraph to the free $\infty$-category it generates has a right adjoint $G:\infty-Cat \rightarrow \Pb$. This right adjoint produces rather complicated polygraphs, indeed if $X$ is an $\infty$-category then $G(X)$, together with the co-unit of adjunction $G(X)^* \rightarrow X$ is constructed inductively as follows: 

\begin{itemize}

\item The $0$-cells of $G(X)$ are exactly the $0$-arrows of $X$.

\item Once $G(X)_n$ and the map $f:G(X)_n^* \rightarrow X$ are constructed. An $(n+1)$-cell of $G(X)$ is given by a triple of: a pair of parallel $n$-arrows $s$,$t$ in $(G(x)_{n})^*$, and an $(n+1)$-arrow $x$ in $X$ between $f(s)$ and $f(t)$. The source and target of the triple $(s,t,x)$ are given by $s$ and $t$, and $f$ is extended by $f(s,t,x)=x$.

\end{itemize}

In particular, there is a terminal object in the category of polygraphs denoted $\Pb 1$, and given by $G(*)$ where $*$ denotes the terminal $\infty$-category (which has only one arrow). And it is a rather rich object: it has one $0$-cell, one $1$-cell, and then one $n$-cell for each pair of parallel $(n-1)$-arrows in $(\Pb 1)^*$. For example $\Pb 1$ has a $2$-cell $A_{n,m}$ from the $n$-times composite to the $m$-times composite of its unique $1$-cell for each pair of integers $n,m \geqslant 0$. This object will play an extremely important role in the theory of good class of polygraphs developed in section \ref{section_goodclass}.

More details about this adjunction, for example its monadicity, can be found in \cite{metayer2016strict}.

}

\block{Finally, a morphism of polygraphs is a monomorphism if and only if it is injective on cells, an epimorphism if and only if it is surjective on cells, and the category of polygraphs admit epi-mono factorization, i.e. every morphism of polygraphs factors uniquely as an epimorphism followed by a monomorphism, and this is achieve by the corresponding factorization on the sets of $n$-cells. This is proved for example in \cite[5.(9)]{makkai2005word}.

The category of polygraphs has all limits and colimits. Colimits are relatively easy to compute due to the following fact: the functor sending a polygraph to the set of its cells (or to the set of its $n$ dimensional cells) preserve colimits.
}

\subsection{Classes of polygraphs}

\block{\label{Def_classofPG}\Def{A \emph{class of polygraphs}, is a full subcategory $J$ of the category of polygraphs which has the following properties:

\begin{itemize}

\item If $X \in J$ and there exists an arrow $f:Y \rightarrow X$ then $Y \in J$.

\item If $Y \in J$ and there exists an epimorphism $Y \twoheadrightarrow X$ then $X \in J$.

\item If $X_i \in J$ for all $i \in I$, then $\coprod_{i \in I} X_i \in J$

\end{itemize}

}

If $J$ is a class of polygraphs, one will sometimes says that $X$ is a ``$J$-polygraph'' for $X \in J$.

It seems that last two conditions of the definition can be replaced by saying that $J$ is stable under all colimits, but this would require a more detailed investigation of limits and colimits in the category of polygraphs that we would prefer to avoid.

}

\block{\label{classofPG_examples}Some examples of this notions:

\begin{itemize}

\item The category $\Pb$ of all polygraphs.

\item The category $\Pb_n$ of polygraphs of dimension at most $n$.

\item The category $\Glob$ of globular polygraphs, i.e. the polygraphs where both the source and target of every cell is a cell.

\item The category $\Ob \Pb$ of ``opetopic polygraphs'' which are those polygraphs for which the target of every cell is again a cell (and not an arbitrary arrow).

\item The category $\Pb^+$ of ``positive polygraphs'' (or ``non-unital'' polygraphs), i.e. the polygraphs such that the source and the target of each cell are non-identity arrows.

\item The categories $\Pb^{s+}$ and $\Pb^{t+}$ respectively of source-positive and target-positive polygraphs which are the polygraphs such that the source (respectively the target) of each cells is a non-identity arrow.

\end{itemize}
}

\block{The intersection of two classes of polygraphs is again a class of polygraphs, for example $\Pb^+ = \Pb^{s+} \cap \Pb^{t+}$. In particular if $J$ is a class of polygraphs one defines $J_n$ as $J \cap \Pb_n$, i.e. the class of polygraphs in $J$ of dimension at most $n$.}

\block{\Th{Every class of polygraphs $J$ admit a terminal object denoted $J1$. Moreover, $J1$ is a sub-polygraph of the terminal polygraph $\Pb 1$, and $J$ is exactly the class of all polygraphs $X$ whose unique map to $\Pb 1$ factors into $J1$. This induces a correspondence between classes of polygraphs and sub-polygraphs of $\Pb 1$.}

\Dem{For any polygraph $V \in J$ the image (in the sense of the epi-mono factorization) of the unique map $V \rightarrow \Pb 1$ is a polygraph in $J$. The set of all sub-polygraphs of $\Pb 1$ which belongs to $J$ is stable under small unions (because $J$ is stable under small co-product and epimorphic image), hence it has a maximal element, which we call $J1$. As $J1$ is an element of $J$, any object whose unique map to $\Pb 1$ factor into $J1$ is in $J$, and conversely, for any object $X \in J$, its image in $\Pb 1$ is a subpolygraph of $\Pb 1$ in $J$ and hence is included in $J 1$, so the unique map from $X$ to $\Pb 1$ factors into $J1$. This shows that $J$ is the class of polygraphs whose map to $\Pb 1$ factors into $J1$. In particular, $J1$ is the terminal object of $J$.

Conversely, if $W \subset \Pb 1$ is any sub-polygraph of $\Pb 1$ then the category of polygraphs over $W$ is a class of polygraphs whose terminal object is $W$, and this proves the correspondence.
}

}

\section{Good class of polygraphs and pasting diagrams}
\label{section_goodclass}

\blockn{The overhaul goal of this section (and in fact of this paper) is to provide tools to show that certain classes of polygraphs are presheaf categories and to use these tools to prove that the category of positive, source-positive or target-positive polygraphs are all presheaf categories.

For example, the category of $1$-polygraphs, i.e. of ordinary oriented graphs is clearly a presheaf category, and more generally, the category $\G$ of globular polygraphs is equivalent to the category of globular sets and is a presheaf category.

It has been known for a long time that the category of $2$-polygraphs is also a presheaf category. This is originally due to S.Schanuel (unpublished) and first appears in \cite{carboni1995connected}. At some point, some people were lead to believe that the category of all polygraphs was itself a presheaf category (this was erroneously claimed in \cite{carboni1995connected} and generalized in \cite{batanin1998computads}) but this was proved to be false in \cite{makkai2008category} and with a more direct argument in \cite{cheng2012direct}. There is also a rather large class of polygraphs which is known to be a presheaf category: the class $\Ob \Pb$ of opetopic polygraphs, or ``many-to-one'' polygraphs is equivalent to the category of opetopic sets, and hence is a presheaf category, but the proof of this is rather indirect\footnote{Since this paper was first written, a more direct proof of this fact have been given in \cite{thanh2018equivalence}}: this follows from \cite{hermida2000weak} which proves that opetopic sets are the same thing as something called multitopic sets and \cite{harnik2008computads} which proves that these multitopic sets are themselves the same thing as many-to-one polygraphs. One could\footnote{One can at least deduce it from proposition 5.4.10 of \cite{henry2018regular}, and the results of the present paper. But we believe it can be proved more directly.} also probably show that Street's orientals (\cite{street1987algebra}) generate a class of polygraphs which is equivalent to the category of semi-simplicial sets, and hence is also a presheaf category.

\bigskip

Also the work of M.Batanin in \cite{batanin2002computads} gives a criterion to show that certain analogue of the notion of polygraphs for more general globular operads are presheaf categories, it applies in particular to the category of polygraphs for Gray categories, or the category of polygraphs for any cellular globular operads. This criterion of M.Batanin is in fact very similar to the one we develop in the present paper, but cannot be used on the new examples of positive polygraphs we study here. The relation between the two criterion is discussed in \ref{Discuss_Rel_Batanin}. Of course, M.Batanin has developed his criterion in a more general framework of globular operads while we work specifically with the operads for strict $\infty$-categories, but, as we will explained at the end of this introduction, we are convinced that the methods used here can be generalized to encompass M.Batanin work, and we hope to come back to this in a future work. }

\blockn{Our general strategy is as follows: One first introduces in section \ref{subsection_Googclass} a notion of ``good class of polygraphs'' (see definition \ref{Def_goodclass}), a good class is essentially a class which is a presheaf category and which admits a good notion of ``pasting diagrams'' such that any arrow in the free $\infty$-category $X^*$, generated by a given polygraph $X$ in that class, can be represented uniquely by a polygraphic map from one of these pasting diagrams to $X$.
>
\bigskip

Following a terminology suggested to us by A.Burroni (that he used in \cite{burroni2012automates}), we will call a representable object of a good class a ``plex'' ( M.Makkai also used the name ``Computope'' in \cite{makkai2005word}) and a pasting diagram will be called a ``polyplex''.

\bigskip

Section \ref{subsec_FRfunctor} is an introduction to the ideas of A.Carboni and P.T.Johnstone from \cite{carboni1995connected} that we will use constantly in the rest of the section. They essentially show how from a good notion of pasting diagram in dimension $n$ one can prove that one has a presheaf category in dimension $n+1$. We illustrate this idea in \ref{Discuss_2PG_epc} on a well known example: there is good notion of pasting diagram for the free category on a graph and that this allows to see that the category of $2$-polygraphs is a presheaf category. The way P.T.Johnstone and A.Carboni arrived from there, in \cite{carboni1995connected}, at a false claim that the category of all polygraphs was a presheaf category was by believing that such a notion of ``pasting diagram'' should exist in all dimensions. But this is false in general, exactly because of the Eckmann-Hilton argument, and it can be an extremely subtle issue in restricted cases.

\bigskip

The pivotal new results of this paper are those of section \ref{subsection_PlexInduction}, and more specifically theorem \ref{Th_mainGoodClass}. More precisely, we take some class of polygraphs $J$ which is good up to dimension $n$. The results of Johnstone and Carboni implies that $J_{n+1}$ is a presheaf category as well. We will define (in \ref{Def_Jnp1Pasting}) a notion of ``$(n+1)$-polyplex'' for this class. Those $(n+1)$-dimensional polyplexes can fails to represent uniquely the arrows of a free $\infty$-category, but our main theorem ( \ref{Th_mainGoodClass}) is that this failure is only up to the automorphisms of the $(n+1)$-polyplexes which preserve the arrow they represents, and $J_{n+1}$ will be a good class in dimension $n+1$ if and only if there is no such automorphisms of its $(n+1)$-polyplexes.

\bigskip

This allows for proofs that certain classes of polygraphs are presheaf categories by induction on the dimension, where at each step one just needs to prove that polyplexes have no automorphisms. That is what we do in subsection \ref{subsection_AutomPD-SPpolyg} to show that the class of source-positive polygraphs is a presheaf category by analyzing more precisely these automorphisms of polyplexes. This also show that target-positive polygraphs form a presheaf category and it implies that positive polygraphs are a presheaf category as well.

\bigskip

Somehow a polyplex can be thought of as encoding a ``generic'' composition operation that one can perform in an $\infty$-category, its underlying polygraph being its ``arity''. These automorphisms of a polyplex preserving the arrow it represents that we mentioned can be thought of as a permutation of the variables of this composition operation that preserves the result, i.e. a form of commutativity of the operation. This is exactly what happens in the situation of the Eckmann-Hilton argument where the composition of two endomorphisms of an identity cell happen to be commutative and this prevents the class of all $2$-polygraphs to be such a good class of polygraphs, and the class of $3$-polygraphs to even be a presheaf category.
}

\blockn{Finally, it is highly probable that under its present form our theorem \ref{Th_mainGoodClass} has no other applications than proving that the class of source/target-positive polygraphs is a good class of polygraphs. But most of its proof does not really use specific properties of the notion of $\infty$-category and we strongly believe that it should have some very important generalizations to similar situations based, for example, on different globular operads and even on different, non-globular, combinatorics. We hope to soon be able to extend this result to, for example, any finitary parametric right adjoint cartesian monad on a presheaf category over a direct indexing category (possibly a directed category having some automorphisms), using the extended notion of polygraphs proposed by M.Shulman\footnote{This idea of M.Shulman unfortunately only appears on the nLab website, on the entry ``computads'', section ``on inverse diagrams'' at \url{https://ncatlab.org/nlab/show/computad#on_inverse_diagrams} from revision 26.}. }

\subsection{Familially representable functors and presheaf categories}
\label{subsec_FRfunctor}

\blockn{In this subsection we recall well known results about familially representable functors, mostly coming from \cite{carboni1995connected}\footnote{Note that this reference contains two important mistakes. They are explained in \cite{carboni2004corrigenda}, and do not affect the part we will use here.},  that will be one of the main tool toward proving that certain classes of polygraphs are presheaf categories. All the proofs of the claims made below can be found in \cite{carboni1995connected}.}

\block{\Def{A functor $F:\Ccal \rightarrow Sets$ is said to be familially representable if up to natural isomorphism it is of the form:

\[ F(X) \simeq  \coprod_{i \in I} Hom(c_i,X), \]

for a set $I$ and a family of objects $c_i \in \Ccal$.
}}

\block{\label{Prop_NTFRfunc}\Prop{The family $(I,(c_i))$ defining a familially representable functor $F$ is uniquely determined (up to unique isomorphism) by the functor $F$. Natural transformations between such functors correspond to morphisms in the category of family $Fam(\Ccal^{op})$, i.e. a morphism $(I,(c_i)) \rightarrow (J,(d_j))$ is a pair of:

\begin{itemize}

\item a function $f:I \rightarrow J$,

\item an $I$-indexed collection of morphisms $f_i:d_{f(i)} \rightarrow c_i$.

\end{itemize}
Composition is simply given by compositions of functions between sets and compositions of the arrows in $\Ccal$.

}
}

\block{\label{Rk:FMrepfct_and_local_initial_obj}
A nice way to see that the family $(I,c_i)$ is well defined from $F$ up to unique isomorphism, is the following observation: A functor $F$ is familially representable if and only if the category of element of $F$, i.e. the category of pairs $X\in \Ccal, v \in F(X)$, is a disjoint union of a set of categories which all have an initial object. The set $I$ is then the set of connected component of this category of elements, and for each such component $i \in I$, its initial object is the pair $(c_i,v_i)$ where $v_i$ is the element of $F(c_i)$ corresponding to the identity of $c_i$ under the isomorphism:

\[ F(c_i) \simeq  \coprod_{j \in I} Hom(c_j,c_i). \]

}

\block{\label{Prop_CartNTFRfunc}\Prop{For a natural transformation $\lambda:F\rightarrow G$ between familially representable functors $F \sim (I,c_i)$ and $G \sim (J,d_j)$ the following conditions are equivalents:

\begin{itemize}

\item The natural transformation $\lambda$ is cartesian, i.e. for any morphism $f:X \rightarrow Y$ in $\Ccal$ the naturality square of $\lambda$:

\[\begin{tikzcd}[ampersand replacement=\&]
F(X) \arrow{d}{\lambda_X} \arrow{r}{F(f)} \& F(Y) \arrow{d}{\lambda_Y} \\
G(X) \arrow{r}{G(f)} \& G(Y) \\ 
\end{tikzcd}
\]

is a pullback square.

\item All the morphisms $\lambda_i:d_{\lambda(j)} \rightarrow c_i$ of the morphism of family corresponding to $\lambda$ are isomorphisms.

\end{itemize}

}}

\block{\label{Prop_ColimOfFRFunc}\Prop{If the category $\Ccal$ has all finite co-limits, then the category of familially representable functors on $\Ccal$ is stable under finite limits in the category of all functors from $\Ccal$ to $Sets$. Moreover in terms of family, these limits can be described as follow:

\begin{itemize}

\item The terminal object is the family with only one object, which is initial in $\Ccal$.

\item If $(I,c_i)$ and $(J,d_j)$ are family their product is given by $(I \times J, c_i \coprod d_j)$.

\item If one has two arrows $(I,c_i) \rightarrow (K,e_k)$ and $(J,d_j) \rightarrow (K,e_k)$ their fiber product is given by:

\[ \left(I \times_K J, c_i \coprod_{e_k} d_j\right) \]

where $I \times_K J$ denotes the ordinary fiber product of sets, i.e. the set of pairs $(i,j) \in I \times J$ which have the same images in $K$, and the letter $k$ just denotes this common image in $K$.

\item If one has two arrows $f,g: (I,c_i) \rightrightarrows (J,d_j)$ their equalizer is given by:

\[ \left( \{ i \in I | f(i)=g(i) \} , coeq(f_i,g_i: d_{f(i)} \rightrightarrows c_i \right)  \]

\end{itemize}

}}

\blockn{The reason familially representable functors are important in order to prove that certain classes of polygraphs are presheaf categories is the following theorem due to A.Carboni and P.Johnstone in \cite{carboni1995connected}  (theorem 4.1 and 4.3 of this reference) which we only state in the restricted case that is of interest to us:  }

\block{\label{Th_JohnstoneCarbonni}\Th{Let $\Ecal$ be a category with a terminal object. Let $F:\Ecal \rightarrow Sets$ be a functor. Then the comma category $Sets/F$ defined below is a presheaf category if and only if $\Ecal$ is a presheaf category and $F$ is familially representable.}

The comma category $Sets/F$ is the category of triples $e \in \Ecal$ , $X$ a set and $f:X \rightarrow F(e)$ where morphisms are the pair of morphisms $e \rightarrow e'$ and $X \rightarrow X'$ which make the natural square commutes.

}

\block{\label{Discus_JohnstoneCarbonni}One can be a little more precise about this theorem: Assume that $\Ecal = Prsh(\Ccal)$ is a presheaf category, and let:

\[ F(X) = \coprod_{i \in I} Hom(\partial a_i ,X)\]

be a familially representable functor from $Prsh(\Ccal)$ to sets, represented by the family of presheaves $(\partial a_i)_{i \in I}$. Then it is not very hard to see that the comma category $Sets/F$ is equivalent to the category of presheaves over the small category $\Ccal^+$ described as follow:

\begin{itemize}

\item $\Ccal^+$ contains $\Ccal$ as a full subcategory and has in addition a set of new objects $(a_i)$ for $i \in I$.

\item The only morphisms in $\Ccal^+$ are the morphisms in $\Ccal$, the identities of the $a_i$, and for each $c \in \Ccal$ and $i \in I$ one has that:

\[ Hom(c,a_i) = \partial a_i (c)\]

composition with morphisms in $\Ccal$ is given by the functoriality of the $\partial a_i$.

\end{itemize}

Moreover a presheaf $P$ on $\Ccal^+$ corresponds to an object $(e,X,f)$ of $Sets/F$ as follows: the object $e \in \Ecal = Prsh(\Ccal)$ is the restriction of $P$ to $\Ccal \subset \Ccal^+$, the sets $X$ is:

\[ X= \coprod_{i \in I} P(a_i)\]

and an element $(i,v \in P(a_i))$ is sent to $F(e)$ by restricting it along the natural map $\partial a_i \rightarrow a_i$ to get a morphism $\partial a_i \rightarrow e$, corresponding to an element of $F(e)$.

Moreover we will only use this direction of the theorem of A.Carboni and P.T.Johnstone. The converse direction (that if $F$ is not familially representable then $Sets/F$ is not a presheaf category) is only important to us at a ``philosophical'' level to show that our approach is essentially optimal. It will not play any concrete role in the present work.
}

\block{\label{Discuss_2PG_epc}As an example, we will use this to show that the category $\Pb_2$ of $2$-polygraphs is a presheaf category following the proof given in \cite{carboni1995connected} (example 4.6). We include this proof only because it is the simplest non-trivial case of the general strategy we will use in this paper and we believe it is an illuminating example.

A $2$-polygraph is given by the data of a $1$-polygraph, i.e. an ordinary graph, together with a collection of $2$-cells, each $2$-cells being attached to a pair of parallel $1$-cells, hence the category of $2$-polygraphs can be described as the comma category $Set/D$ where $D$ is the functor from the category of graphs to the category of sets, which maps a graph $G$ to:

\[ D(G)= \left\lbrace \text{Pair of Parallel arrows in the free category $G^*$ of path in $G$} \right\rbrace \]

Now, the category $\Pb_1$ of $1$-polygraphs (i.e. of graphs) is equivalent to a presheaf category: the category with two objects $\underline{P}_0$ and $\underline{P}_1$ and with only two non-identity arrows $s,t:\underline{P}_0 \rightrightarrows \underline{P}_1$. So in order to conclude we need to prove that the functor $D$ above is familially representable.

One starts with the simpler functor $P$ which maps a graph $G$ to the set of all arrows of free category $G^*$ of paths in $G$. It is familially representable, indeed the set of path of length $n$ is the same as $Hom(\underline{P}_n,G)$, with $\underline{P}_n$ being the graph:

\[\underline{P}_n :  a_0 \rightarrow a_1 \rightarrow \dots \rightarrow a_n \]

hence: 

\[ P(G)= \coprod_{n=0}^{\infty} Hom( \underline{P}_n,G) \]

Similarly, the functor $D$ is familially represented by the following objects:

 \[ \begin{tikzcd}[ampersand replacement=\&]
\& \& a^s_1 \arrow{r} \& \dots \arrow{r} \& a^s_{n-1} \arrow{dr} \&  \\
\partial \underline{A}_{n,m} : \& a_0 \arrow{ur} \arrow{dr} \& \& \& \& a^s_n=a^t_m \\
\& \& a^t_1 \arrow{r}\& \dots \arrow{r} \&  a^t_{m-1} \arrow{ur} \&   \\
\end{tikzcd}
\]

each object $\partial \underline{A}_{n,m}$ representing the functor sending a graph $G$ to the set of pairs of parallel arrows $(s,t)$ in $G^*$ with $s$ of length $n$ and $t$ of length $m$.

This concludes the proof that the category of $2$-polygraphs is a presheaf category, one moreover gets explicitly the category on which they are presheaves, its objects are:

\begin{itemize}

\item  $\underline{P}_0$ the polygraph with only one $0$-cell: $*$.
\item  $\underline{P}_1$ the polygraph with two $0$-cells and one $1$-cell between them:

\[ \bullet \rightarrow \bullet \]

\item for all integer $n,m \geqslant 0$,  $\underline{A}_{n,m}$ the polygraph

 \[ \begin{tikzcd}[ampersand replacement=\&]
 \& a^s_1 \arrow{r} \& \dots \arrow[Rightarrow,shorten >= 10pt,shorten <= 10pt]{dd}{A_{n,m}} \arrow{r} \& a^s_{n-1} \arrow{dr} \&  \\
 a_0 \arrow{ur} \arrow{dr} \& \& t \& \& a^s_n=a^t_m \\
 \& a^t_1 \arrow{r}\& \dots \arrow{r} \&  a^t_{m-1} \arrow{ur} \&   \\
\end{tikzcd}
\]

\end{itemize}

The morphisms being the morphisms of polygraphs between them, but more explicitly there is no non-identity morphisms between the $\underline{A}_{n,m}$ and the morphisms from $\underline{P}_0$ and $\underline{P}_1$ to $\underline{A}_{n,m}$ are just the graph morphisms to $\partial \underline{A}_{n,m}$.

}

\subsection{Good classes of polygraphs}
\label{subsection_Googclass}

\block{As pointed out by M.Makkai in \cite{makkai2005word}, we do not just want a class of polygraphs $J$ to be a presheaf category, we want it to be what he calls an ``effective presheaf category''. An effective category is a category $\Ccal$ endowed with a functor $F : \Ccal \rightarrow Set$ called the forgetful functor. It is common to ask additional conditions on $F$, like being faithful or an isofibration, but those assumptions plays no role here (although they will be satisfied on all the examples we will encounter). Two effective categories are said to be equivalent if there is an equivalence between them which commutes to the forgetful functor up to natural isomorphism.

The category of polygraphs is an effective category with the forgetful functor mapping a polygraph $X$ to the set of all its cells. Presheaf categories are effective categories under the forgetful functor:

\[ \Fcal \mapsto \coprod_{c \in \Ccal} \Fcal(c) \]

\Def{One says that a class of polygraphs is an effective presheaf category if it is equivalent, as an effective category, to a presheaf category.}

More explicitly this means that there is a full subcategory $I$ of the class $J$ such that the natural functor from $J$ to $Prsh(I)$ is an equivalence and such that the functor from $J$ to sets sending a polygraph to the set of its cell is familially representable by the family $I$:

\[ \{ \text{Cells of } X\} \simeq \coprod_{ \underline{i} \in I} X( \underline{i}) \]

For example, $\Pb_1$ is obviously an effective presheaf category, and the discussion in \ref{Discuss_2PG_epc} shows that $\Pb_2$ is an effective presheaf category.

}

\block{\Def{If a class of polygraphs $J$ is an effective presheaf category, a representable object is called a $J$-plex.
}

Plexes (and latter polyplexes) will always be denoted with an underlined letter (like the $\underline{i}$ in the previous paragraph).

As mentioned above, this terminology is due to A.Burroni in \cite{burroni2012automates}. They have also been called ``computopes'' by M.Makkai in \cite{makkai2005word} who defined them for a general class of polygraphs.

The proposition below will show that to each cell of the terminal $J$-polygraph $J1$ corresponds a $J$-plex which is unique up to unique isomorphism, and that $J$-plexes are classified by the cells of $J1$.

}

\block{\label{prop_effPreshv}\Prop{If a class of polygraphs $J$ is an effective presheaf category, with $I$ the category of $J$-plexes, i.e. $J \simeq Prsh(I)$ as effective categories then:

\begin{enumerate}

\item There is a bijection between the $J$-plexes and the cells of $J1$. 

\item For each cell $c \in J1$ the corresponding $J$-plex $\underline{x} \in I$ is characterized by the following universal property:

\[ Hom(\underline{x},X) = \{ a \in X | \pi(a)=c \text{ where $\pi$ is the unique map $ \pi:X \rightarrow J1$ }\}\]

\item Any $J$-plex $\underline{x}$ has only one cell of dimension $n$, denoted $x$, no cells of dimension higher than $n$, and a finite number of cells in total. Moreover $\underline{x}$ is generated by $x$ in the sense that there is no strict sub-polygraphs of $\underline{x}$ containing $x$.

\item The category $I$ of $J$-plexes is a directed category: every non-identity morphism goes from an object of dimension $n$ to an object of dimension strictly higher. Moreover all the slices of $I$ are finite directed categories.

\item Any subclass of polygraphs of $J$ is also an effective presheaf category. Moreover if $J' \subset J$ is a subclass then the $J'$-plexes are exactly the $J$-plexes that belongs to $J'$, or equivalently, the $J$-plexes corresponding to cells of $J'1$.

\end{enumerate}

}

The cell $x$ of $\underline{x}$ is called the universal cell of $x$, or sometimes the top cell as it is the only cell of maximal dimension.

\Dem{\begin{enumerate}

\item One has that for each polygraph $X$, the set of cells of $X$ is (isomorphic) to:

\[ \coprod_{\underline{i} \in I} Hom(\underline{i},X) \]

In particular, taking $X= J1$, one has that for all $\underline{i} \in I$ $Hom(\underline{i},J1)=\{*\}$ as $J1$ is terminal in $J$. And hence the set of cells of $J1$ can be written as:

\[ \coprod_{\underline{i} \in I} \{*\} \simeq I \]

\item For a general polygraph $X$ in $J$, the unique map from $X$ to $J1$ identifies it at the level of the sets of all cells:

\[ \coprod_{\underline{i} \in I} Hom(\underline{i},X) \rightarrow  \coprod_{\underline{i} \in I} \{*\} \simeq I \]

hence $Hom( \underline{i},X)$ identifies with the set of cells of $X$ which are sent to the cell of $J1$ corresponding to $\underline{i}$.

\item The identity map of a $J$-plex $\underline{x}$ corresponds (by its universal property) to a cell $x$ of $\underline{x}$ which is mapped by the map $ \underline{x} \rightarrow J1$ to the cell $c$ corresponding to $\underline{x}$. Hence $x$ is $n$-dimensional if $c$ has dimension $n$.

Let $V$ be any sub-polygraph of $\underline{x}$ which contains $x$, and let $i:V \rightarrow \underline{x}$ be the inclusion map, because $V$ contains $x$ there is a map $p$ from $\underline{x}$ to $V$ corresponding to this cell by the universal property of $\underline{x}$, but the functoriality of the universal property gives immediately that $i \circ j =Id_{\underline{x}}$ hence as $i$ is an inclusion this implies that $V=X$.

This proves that $\underline{x}$ is generated by $x$ in the sense that every sub-polygraphs that contains $x$ is equal to $\underline{x}$.

One proves by induction on arrows that in any polygraph any $n$-arrow $f$ is contained in a finite sub-polygraph of dimension at most $n$ and whose only $n$-cells are those appearing in $f$. Applying the first part of this point to such a subpolygraph containing $x$ implies that $\underline{x}$ is finite, of dimension $n$ and $x$ is the unique $n$-cell of $\underline{x}$. 

\item Morphisms from an arbitrary object $\underline{y} \in I$ to a fixed object $\underline{x} \in I$ are in one to one correspondence with the cells of $\underline{x}$. So for all of them except one, $\underline{y}$ have dimension strictly smaller than the dimension of $\underline{x}$, and the last one corresponds to the cell $x$, i.e. the identity of $\underline{x}$. This proves that the category of plexes is directed (by dimension). Moreover as morphisms $\underline{y} \rightarrow \underline{x}$ in the category of plexes corresponds to cells of $\underline{x}$, this shows that the slice category at $\underline{x}$ is finite.

\item Let $J' \subset J$ is a sub-class of polygraphs. For each cell $c \in J'1 \subset J1$, there is a unique $J$-plex with a morphism $\underline{x} \rightarrow J'1 \subset J1$ sending $x$ to $c$. This plex $\underline{x}$ is in $J'$ as its unique map to $J1$ factor into $J'1$, and it is the $J$-plex corresponding to the cell $c \in J' 1 \subset J1$. So for a $J$-plex it is equivalent to be in $J'$ and that the corresponding cell of $J1$ is in $J'1$. Let $I'$ be the full subcategory of $I$ of such plexes.

As a presheaf over $I$, the object $J' 1$ hence corresponds to the presheaf sending a plex $\underline{x} \in I$ to $\{*\}$ if $x \in I'$ and $\emptyset$ otherwise (which is indeed a presheaf, i.e. $I'$ is downward closed in $I$). $J'$ corresponds to the slice of $J$ over this object, but in terms of presheaves those are exactly presheaves over $I'$, extended by $\emptyset$ on the object of $I$ not in $I'$.

\end{enumerate}
}
}

\block{\label{Def_goodclass}As we mentioned in \ref{Th_JohnstoneCarbonni}, a key property for proving that some class of polygraphs is a presheaf category is that the free $\infty$-category functor should be familially representable on this class. For this reason we will introduce the following notion of ``good class of polygraphs'', which will be our typical induction hypothesis in a proof that a given class of polygraphs is a presheaf category.

\Def{A class of polygraphs $J$ is said to be a \emph{good} class of polygraphs if:

\begin{itemize}

\item $J$ is an effective presheaf category.

\item The functor $X \mapsto X^*$ from $J$ to the category of Sets, which maps a polygraph $X$ to the set of all arrows of the free $\infty$-category $X^*$ is familially representable.

\end{itemize}

}

}

\block{\label{DEF_polyplex}\Def{If $J$ is a good class of polygraphs, the objects that familially represent $X \mapsto X^*$, i.e. the $\underline{p}$ such that:

\[ X^* \simeq \coprod_{\underline{p}} Hom(\underline{p},X) \]

Are called the $J$-polyplexes (or $J$-pasting diagrams).

}

Here again, the terminology ``polyplexes'' comes from A.Burroni in \cite{burroni2012automates}. Similarly to plexes, polyplexes will also be always denoted with an underlined letter.

}

\block{\label{Prop_goodclass}\Prop{If $J$ is a good class of polygraphs then:

\begin{itemize}

\item The $J$-polyplexes are in bijection with the arrows of $(J1)^*$.

\item The $J$-polyplex $\underline{p}$ corresponding to an arrow $c \in (J1)^*$ is characterized by the universal property:

\[ Hom(\underline{p},X) = \{ x \in X^* | \text{s.t. } \pi^*(x)=c\} \]

where $\pi^*$ denotes the unique polygraphic map $\pi^*:X^* \rightarrow (J1)^*$.

\item If $\underline{p}$ is a polyplex, and $p \in \underline{p}^*$ is the arrow corresponding to the identity map of $\underline{p}$ in the isomorphisms above, then $\underline{p}$ is generated by $p$ in the sense that there is no proper subpolygraphs $Y \subset \underline{p}$ such $p \in Y^*$. Moreover $\underline{p}$ is finite and of dimension the dimension of $p$.

\item If $J' \subset J $ is any subclass of polygraphs, then $J'$ is a good class of polygraphs and the $J'$-polyplexes are exactly the $J$-polyplexes that belong to $J'$, or equivalently, the $J$-polyplexes corresponding to arrows of $(J'1)^* \subset (J1)^*$.

\end{itemize}

}

The arrow $p \in \underline{p}^*$ is called the universal arrow, and we will always denote it by the same letter as the polyplex. The universal property of polyplexes can be rephrased as: 

If $X$ is a $J$-polygraph and $a \in X^*$ then there is a unique polyplex $\underline{p}$ and map $\chi_a: \underline{p} \rightarrow X$ such that $\chi_a^* p = a$.

\bigskip

Note that $p \in \underline{p}^*$ is part of the structure of ``being a polyplex'' in the sense that it is the choice of this arrow that explains how the functor $Hom(\underline{f}, \_)$ appears as a subfunctor of $( \_)^*$. In fact, in \ref{ExIsomorphicPolyplex} and \ref{Ex:SimonForest}  we will construct two pairs of examples of polyplexes in the (good) class of positive $3$-polygraphs which are isomorphic as polygraphs but with different universal arrows, which makes them non-isomorphic as polyplexes.

In terms of the discussion in \ref{Rk:FMrepfct_and_local_initial_obj}, polyplexes are the pairs $(\underline{f},f)$ which are initial in their connected component of the category of pairs $(X \in J, a \in X^*)$.

Finally comparing the second point of proposition \ref{prop_effPreshv} and \ref{Prop_goodclass} immediately show that plexes are exactly the polyplexes corresponding to generators of $(J1)^*$, and the ``universal arrow'' of a plex (when seen as a polyplex) is its universal cell. Plexes are also exactly the polyplexes whose universal arrows is a generator.

\Dem{This is very similar to the proof of proposition \ref{prop_effPreshv}:

By assumption, the arrows of $(J1)^*$ are in bijection with the coproduct over the $J$-polyplexes of the $Hom(\underline{f},J1)=\{*\}$, hence this defines a bijection between arrows of $(J1)^*$ and the set of $J$-polyplexes. Moreover for any object $X$ one has $\pi^* : X^* \rightarrow J1^*$ which corresponds to 

\[\coprod_{\underline{v} \text{ a polyplex}} Hom(\underline{v},X) \rightarrow \coprod_{\underline{v} \text{ a polyplex}} Hom(\underline{v},J1) \simeq (J1)^* \]

which immediately gives the universal property of $\underline{v}$ claimed in the proposition.
If $Y \subset \underline{p}$ is a subpolygraph such that $p \in Y^*$, then the universal property of $\underline{p}$ implies that the identity map of $\underline{p}$ factors into $Y$, which implies that $Y = \underline{p}$. As any arrow of dimension $n$ in a polygraph is contained in a finite subpolygraphs with only cells of dimension $\leqslant n$ (easy by induction on arrows) this implies that $\underline{p}$ is finite and of dimension $n$.

Finally, if $J' \subset J$ is any subclass of polygraphs of $J$, then $J'$ is an effective presheaf category by \ref{prop_effPreshv}, and for any $X \in J'$ one has that:

\[ X^* = \coprod_{\underline{v} \text{ a $J$-polyplex}} Hom(\underline{v},X) \]

but any $J$-polyplex admitting a morphism to $X$ is in $J'$ so one has that:

\[X^* = \coprod_{\underline{v} \text{ a $J'$-polyplex}} Hom(\underline{v},X) \]

where the $J'$-polyplexes are defined as being the $J$-polyplexes that belongs to $J'$.

}
}

\block{In all this paper we will identify $J$-plexes with cells of $J1$ and $J$-polyplexes with arrows of $(J1)^*$. When doing so, cells and arrows of $(J1)^*$ will also be denoted with an underlined letter as well. The results of the present section show that these conventions are compatible to essentially everything the author has been able to think of:

\begin{itemize}

\item If $\underline{a}$ is a generator of $(J1)^*$ then the corresponding plex and the corresponding polyplex are the same.

\item The dimension of an arrow (or a generator) of $(J1)^*$ is the same as the dimension of the underlying polygraph of the corresponding polyplex (or plex).

\item If $J' \subset J$ is a subclass of polygraphs, then a $J$-plex (resp. polyplex) $\underline{p}$ is in $J'$ if and only if the corresponding cell (resp. arrow) is in $J'1$ (resp. $(J'1)^*$).

\end{itemize}

}

\block{The functor $X \mapsto X^*$ does not just take values in the category of sets, but in the category of $\infty$-categories. In particular this functor is endowed with all the structure of an $\infty$-category. Because of proposition \ref{Prop_NTFRfunc}, a structure on a familially representable functor will manifest itself as operations on the representing family, i.e. on the $J$-polyplexes. At the level of indexing set of the family, those operations correspond to the structure of $\infty$-category on $(J1)^*$ under the identification of polyplexes with arrows of $(J1)^*$, but one also get morphisms between the polyplexes:

\bigskip

The source and targets maps $\pi^{\epsilon}_k : X^* \rightarrow X^*$ are represented by morphisms of families, i.e. for each polyplex $\underline{p}$ one has polyplexes $\pi^{\epsilon}_k(\underline{p})$ and morphisms from $\pi^{\epsilon}_k(\underline{p}) \rightarrow \underline{p}$, such that if an arrow of $X^*$ is represented by a morphism $\underline{p} \rightarrow X$ then its $k$-source or $k$-target are represented by the maps $\pi^{\epsilon}_k \underline{p} \rightarrow \underline{p} \rightarrow X$.

\bigskip

Similarly, the composition operation $\_ \#_k \_$ can be seen as an operation $X^* \times_{X^*} X^* \rightarrow X^*$ where the fiber product is over the two maps $\pi^-_k$ and $\pi^+_k$ from $X^*$ to $X^*$. Because of proposition \ref{Prop_ColimOfFRFunc}, this fiber product is a familially representable functor, represented by the family of $\underline{f} \coprod_{\underline{p}} \underline{g}$ where $\underline{f},\underline{g}$ and $\underline{p}$ are polyplexes such that $\pi^+_k \underline{f} = \pi^-_k \underline{g}=\underline{p}$ and the coproduct is taken along the maps $\underline{p}=\pi^+_k \underline{f} \rightarrow  \underline{f}$ and $\underline{p}= \pi^-_k \underline{g} \rightarrow \underline{g}$ mentioned above. And hence the $k$-composition operation on arrows of $X^*$ is encoded by maps $\underline{f} \#_k \underline{g} \rightarrow \underline{f} \coprod_{\underline{p}} \underline{g}$ where, following our convention, $\underline{f} \#_k \underline{g}$ is just the polyplex obtained from the composition operations on arrows of $(J1)^*$. These operations satisfies some associativity, exchange and compatibility conditions just translating the corresponding axioms on $X^*$ which we will not list.

\bigskip

One should be aware that there is a small collision in our notations: If $\underline{p}$ is a polyplex with $p$ its universal arrow, then $\pi^{\epsilon}_k p$ can potentially mean both the arrow $\pi^{\epsilon}_k p \in \underline{p}^*$ or the universal arrow of the polyplex $\pi^{\epsilon}_k \underline{p}$. One should however note that the former is the image of latter by the natural map $\pi^{\epsilon}_k \underline{p} \rightarrow \underline{p}$ so this collision is not completely a bad thing, though it might be good to occasionally precise if we are talking of the map  $\pi^{\epsilon}_k p \in \underline{p}^*$ or $\pi^{\epsilon}_k p \in \pi^{\epsilon}_k \underline{p}^*$. Especially that there are some cases where the map $\pi^{\epsilon}_k \underline{p} \rightarrow \underline{p}$ is not a monomorphism (see the examples in \ref{ExIsomorphicPolyplex}). A similar collision happen with $f \#_k g \in \left( \underline{f} \coprod_{\underline{p}} \underline{g} \right)^*$ and $f \#_k g \in \left( \underline{f} \#_k \underline{g} \right)^*$, with the former being the image of the latter by the comparison map: $\underline{f} \#_k \underline{g} \rightarrow \underline{f} \coprod_{\underline{p}} \underline{g}$. In this case, as explained in the paragraph, it will be shown that this comparison map is always an isomorphism.
}

\block{\label{Def_CondC} We will prove at some point latter (corollary \ref{Cor_CondCforGoodClass}) that in a good class of polygraphs, the maps $\underline{f} \#_k \underline{g} \rightarrow \underline{f} \coprod_{\underline{p}} \underline{g}$ mentioned above are always isomorphisms, i.e. (because of proposition \ref{Prop_CartNTFRfunc}) that all the composition operations:

\[ \_ \#_k \_ : X^* \times_{X^*} X^* \rightarrow X^* \]

are cartesian natural transformations in $X$. In the meantime we will say that a class of polygraphs $J$ satisfies condition $(C)$ if the composition operations are cartesian, i.e. if all the maps $\underline{f} \#_k \underline{g} \rightarrow \underline{f} \coprod_{\underline{p}} \underline{g}$ mentioned above are isomorphisms.
}

\subsection{Plexes and polyplexes in dimension $n+1$}
\label{subsection_PlexInduction}

\block{\label{Def_AutPD}In this subsection, we consider a class of polygraphs $J$ such that $J_n$ is a good class of polygraphs, satisfying condition $(C)$ of \ref{Def_CondC}. Using the theorem of P.T.Johnstone and A.Carboni that we quoted in \ref{Th_JohnstoneCarbonni}, we will deduce that $J_{n+1}$ is an effective presheaf category (proposition \ref{Prop_GoodnImpPrenpp}). We will then define a notion of $J_{n+1}$-polyplexes (this will be done in \ref{Def_Jnp1Pasting}). Informally,  $J_{n+1}$-polyplexes are going to be defined as formal compositions (gluing) of $J_{n+1}$-plexes using the composition operations of $\infty$-categories.

It will be clear from definition \ref{Def_Jnp1Pasting} that the $J_{n+1}$-polyplexes of dimension smaller than $n$ are just the $J_n$-polyplexes in the sense of definition \ref{DEF_polyplex}, and at the end of this subsection it will appears that when $J_{n+1}$ is indeed a good class of polygraphs the $J_{n+1}$-polyplexes of definition \ref{Def_Jnp1Pasting} are the same as the polyplexes in the sense of definition \ref{DEF_polyplex}. In the meantime, as we do not know yet whether or not $J_{n+1}$ is a good class of polygraphs, whenever we speak of $J_{n+1}$-polyplexes we are always referring to the definition of the present section (\ref{Def_Jnp1Pasting}), which as far as we know is the only meaningful one, and for $J_n$-polyplexes the two definitions are already known to be equivalent.

Exactly as the polyplexes of definition \ref{DEF_polyplex}, a $J_{n+1}$-polyplex $\underline{v}$ will actually be a pair $(\underline{v},v)$, with $\underline{v}$ a polygraph (the ``underlying polygraph'') and $v$ an arrow of $\underline{v}^*$, one still calls $v$ the ``universal arrow'' of $\underline{v}$. In particular:

\Def{An automorphism of a $J_{n+1}$-polyplex $\underline{v}$ is an automorphism $w$ of the underlying polygraph of $\underline{v}$ such that $w^* v = v$. One denotes by $G_{\underline{v}}$ the group of automorphisms of a polyplex $(\underline{v},v)$.}

Note that if these were the polyplexes of a good class of polygraphs (in the sense of definition \ref{DEF_polyplex}) one will always have $G_{\underline{v}}=\{1\}$ because of the universal property of polyplexes. But this will not necessarily be the case for our $J_{n+1}$-polyplexes, in fact the group $G_{\underline{v}}$ represents exactly the obstruction for $\underline{v}$ to be a polyplex in the sense of definition \ref{DEF_polyplex}:

The main result of this section (theorem \ref{Th_mainGoodClass}) will be that under the assumption that $J_n$ is a good class of polygraphs satisfying condition $(C)$, one has for any $J_{n+1}$-polygraphs $X$ that:

\[ X^* = \coprod_{\underline{v} \text{ a $J_{n+1}$-polyplex}} Hom(\underline{v},X)/G_{\underline{v}} \]

To put it another way: for each cell $f \in X^*$ there is a unique $J_{n+1}$-polyplex $(\underline{p},p)$ and a map $\chi_f: \underline{p} \rightarrow X$ such that $f=\chi_f^* p$, but with $\chi_f$ being only unique up to an element of $G_{\underline{p}}$.  

In particular $J_{n+1}$ will be a good class of polygraphs if and only if all the automorphism groups $G_{\underline{v}}$ are trivial. Moreover, it will also be the case that when this happens, condition $(C)$ automatically holds for $J_{n+1}$, and this will prove by induction that condition $(C)$ in fact holds for all good class of polygraphs (Corollary \ref{Cor_CondCforGoodClass}).

In the next subsection we will study in more details the properties of automorphisms of polyplexes and we will show in particular that they are indeed trivial for the class of source-positive polygraphs.
}

\block{For the rest of this section one will typically assume that $J$ is a class of polygraphs such that $J_n$ is a good class satisfying condition $(C)$ of \ref{Def_CondC}. }

\block{We first prove that $J_{n+1}$ is an effective presheaf category:

\label{Prop_GoodnImpPrenpp}\Prop{Let $J$ be a class of polygraphs such that $J_{n}$ is a good class of polygraphs, then $J_{n+1}$ is an effective presheaf category.

Moreover the $J_{n+1}$-plexes are exactly:
\begin{itemize}

\item The $J_{n}$-plexes.

\item The pairs of parallel $J_n$-polyplexes $(\underline{a},a)$ and $(\underline{b},b)$ glued together along the identifications $\pi^-_{n-1} \underline{a} = \pi^{-}_{n-1} \underline{b}$ and $\pi^+_{n-1} \underline{a} = \pi^{+}_{n-1} \underline{b}$ with just an additional $(n+1)$-cell between $a$ to $b$, such that the resulting polygraph is in $J_{n+1}$.

\end{itemize}
}

\Dem{Let $T$ be the class of all $(n+1)$-polygraphs whose underlying $n$-polygraph belongs to $J_n$. Then $J_{n+1} \subset T$ so that it is enough to prove that $T$ is an effective presheaf category because of \ref{prop_effPreshv}. We will use the theorem of A.Carboni and P.T.Johnstone that we quoted in \ref{Th_JohnstoneCarbonni}, indeed $T$ can be described as the category $Sets/F$ where $F$ is the functor from $J_n$ to sets which maps each polygraph $X$ to the set of pairs of parallel $n$-arrows of $X^*$, this can be described as:

\[ F(X):= \left\lbrace (x,y) \in (X^*)^2  \left| \begin{array}{ c } \pi^-_{n}(x)=x \text{ ; } \pi^-_{n}(y)=y  \\ \forall \epsilon, \pi^{\epsilon}_{n-1}(x)=\pi^{\epsilon}_{n-1}(y) \end{array} \right.  \right\rbrace \]

As all the $\pi^{\epsilon}_k$ are natural transformations of the functor $X \mapsto X^*$, this functor $F$ is familially representable as a finite limit of familially representable functors (see \ref{Prop_ColimOfFRFunc}). Using the description of finite limits of familially representable functors given in \ref{Prop_ColimOfFRFunc}, the family representing $F$ is exactly given by the set of pairs of parallel $J_n$-polyplexes $\underline{a},\underline{b}$, with as objects the gluing of $\underline{a}$ and $\underline{b}$ along the identifications $\pi^{\epsilon}_{n-1} \underline{a} = \pi^{\epsilon}_{n-1} \underline{b}$.

Hence ones immediately deduces that $T$ is a presheaf category, we need to show that it is an effective presheaf category. This follows from the discussion in \ref{Discus_JohnstoneCarbonni}:  If $X$ is a $T$-polygraph represented as an object of $Sets/F$ with $X_n \in J_n$ and $S \rightarrow F(X)$ being the set of $(n+1)$-cells of $X$, the total set of cells is $S \coprod |X_n|$ (where $|X_n|$ denotes the set of cells of $X_n$). The category that we get from \ref{Discus_JohnstoneCarbonni} has for objects the $J_n$-plexes and one additional object for each pair of polyplexes that appears in the representation of $F$. Then $|X_n|$ corresponds to the coproduct over all the $J_n$-plexes of $Hom(\underline{v},X)$ because $J_n$ is an effective presheaf category and $S$ corresponds to the coproduct over the new objects $\underline{w}$ of $Hom(\underline{w},X)$ by \ref{Discus_JohnstoneCarbonni}.

Moreover the $T$-plexes are exactly the objects described in the proposition: the $J_n$-plexes and one additional object for each pair $(\underline{u},\underline{v})$ of parallel $J_n$-polyplexes which represent exactly the $(n+1)$-cell of $T1$ between the arrows $\underline{u}$ and $\underline{v}$ hence is exactly as described in the proposition. Finally, the $J_{n+1}$-plexes are the $T$-plexes which belong to $J_{n+1}$ by proposition \ref{prop_effPreshv}.
}
}

\block{We still assume that $J_n$ is a good class of polygraphs satisfying condition $(C)$. The general idea of our definition of $J_{n+1}$-polyplexes is quite simple: one wants them to be defined inductively by the fact that each $J_{n+1}$-plex is a $J_{n+1}$-polyplex and that if two $J_{n+1}$-polyplexes are ``composable'' then their composite (constructed as a gluing in the spirit of condition $(C)$ of \ref{Def_CondC}) should also be a $J_{n+1}$-polyplex. This will only be possible because the boundary of a $J_{n+1}$-polyplex will be made of two $J_n$-polyplexes, which have no automorphisms and so when we compose two $J_{n+1}$-polyplexes together, there is no ambiguity on how we glue their boundary together. Without this assumption there would be several ways to form the same composition $\underline{x} \#_n \underline{y}$ by twisting it by an automorphism the identification $\pi^+_n \underline{x} \simeq \pi^-_n \underline{y}$, and one could end up with possibly more isomorphism classes of polyplexes than arrows of $(J_{n+1} 1)^*$.

So, as one wants to have one $J_{n+1}$-polyplex for each arrow of $(J_{n+1} 1 )^*$, and one wants them to be constructed by induction on arrows, the simplest way is to construct them as the image of a certain morphism of $\infty$-categories $F:(J_{n+1} 1 )^* \rightarrow \Dcal$ where $\Dcal$ is an $\infty$-category (in fact an $n+1$-category). The objects of $\Dcal$ will be ``pre-polyplexes'', i.e. polygraphs with a marked $(n+1)$-arrow which are composed by gluing, and $F$ will just be defined as sending the generators of $(J_{n+1} 1)^*$ to the corresponding $J_{n+1}$-plexes. So the main point is to construct this $\infty$-category $\Dcal$.

For technical reason, we will need an $\infty$-category $\Dcal_X$ depending functorially on a polygraph $X \in J_{n+1} $ whose $(n+1)$-arrows are ``pre-polyplexes'' with a map to $X$, and which will comes with a natural morphism $X^* \rightarrow \Dcal_X$. The $\Dcal$ we mentioned above corresponds to the case $X=J_{n+1} 1$, but this more general $\Dcal_X$ will be convenient to show that any arrow of $X^*$ can be represented by a $J_{n+1}$-polyplex. }

\block{\label{DefOfDcalX} \textbf{Definition of $\Dcal_X$:}

The arrows of $\Dcal_X$ are equivalence classes of triples $(v,\alpha,\lambda)$ where:

\begin{itemize}
\item $v$ is a finite polygraph in $J_{n+1}$,

\item $\alpha$ is an arrow of $v^*$ of dimension at most $n+1$,

\item $\lambda$ is a polygraphic morphism $\lambda: v \rightarrow X$.

\end{itemize} 

Two such triples $(v,\alpha,\lambda)$ and $(v',\alpha',\lambda')$ are equivalent if there exists an isomorphism $\theta:v \rightarrow v'$ such that $\theta^*(\alpha) = \alpha'$ and $\lambda' \theta = \lambda$.

\bigskip

Sources and targets in $\Dcal_X$ are defined as follows: If $(v,\alpha,\lambda)$ is an arrow of $\Dcal_X$ and $i \leqslant n$ then the arrow $\pi^{\epsilon}_i(\alpha)$ is of dimension at most $n$, so it belongs to the $n$-skeleton of $v$ which is an object of $J_n$, hence there is a uniquely defined $J_n$-polyplex $\underline{x}$ and a map $\chi_{\pi^{\epsilon}_i(\alpha)}: \underline{x} \rightarrow v$ corresponding to this arrow.

One defines:  $\pi^{\epsilon}_i(v,\alpha,\lambda) = ( \underline{x}, x, \lambda \circ \chi_{\pi^{\epsilon}_i(\alpha)})$.

For any arrow $x$ of $\Dcal_X$ one defines $\pi^{\epsilon}_{n+1}(x)=x$ (i.e. $\Dcal_X$ will be an $(n+1)$-category).

\bigskip

All compositions $x \#_i y$ for $i>n$ are trivial as $\Dcal_X$ is an $(n+1)$-category. Compositions for $i \leqslant n$ are defined as follow:

If $\pi^+_i(v,\alpha,\lambda) = \pi^-_i(v',\alpha',\lambda')$ in $\Dcal_X$, it means that there is an isomorphism between the $J_n$-polyplexes representing $\pi^+_i(\alpha)$ and $\pi^-_i(\alpha')$ which send the universal arrow to the universal arrow, but such an isomorphism, when it exists, is unique, and is in fact equal to the identity as $J_n$ is a good class of polygraphs and those are $J_n$-polyplexes, hence one can define:

\[ (v,\alpha,\lambda) \#_i (v',\alpha',\lambda') := \left( v \coprod_{\underline{x}} v', \alpha \#_i \alpha', (\lambda,\lambda') \right) \]

where $\underline{x}$, with its map to $v$ and $v'$ is the $J_n$-polyplex representing $\pi^+_i \alpha$ and $\pi^-_i \alpha'$. The composition $\alpha \#_i \alpha'$ makes sense because $\alpha$ and $\alpha'$ are arrows of $v^*$ and $v'^*$ respectively hence they both belong to $\left(  v \coprod_{\underline{x}} v'\right)^*$, and as $\pi^+_i(\alpha)$ and $\pi^-_i(\alpha)$ both are the image of $x \in \underline{x}^*$, they are equal in the pushout. Similarly $\lambda$ and $\lambda'$ coincide on $\underline{x}$ because of the equality $\pi^+_i(v,\alpha,\lambda) = \pi^-_i(v',\alpha',\lambda')$ on the third component.

\Prop{$\Dcal_X$ defined above is an $(n+1)$-category. Its underlying $n$-category is isomorphic to the underlying $n$-category of $X^*$.}

\Dem{If we restrict ourselves to the arrows of $\Dcal_X$ which are of the form $(\underline{x},x,v:\underline{x}\rightarrow X)$ for $\underline{x}$ a $J_n$-polyplex one gets exactly the $k$-arrows of $X^*$ for $k\leqslant n$, with the correct sources and targets, and the compositions law is correct exactly because we assumed that $J_n$ satisfies condition $(C)$ of \ref{Def_CondC}. Moreover the $\pi^{\epsilon}_k$ for $k \leqslant n$ takes values in this class of arrows, so the ``underlying $n$-category'' of $\Dcal_X$ is indeed an $\infty$-category simply because it is isomorphic to the underlying $n$-category of $X^*$.

Moreover, as all the other $\pi^{\epsilon}_k$ (for $k>n$) are the identity on $\Dcal_X$ all the ``globular relations'' between the $\pi^{\epsilon}_k$  are automatically deduced from the globular relations in $X^*$. At this point we only have to check that the composition operations on $\Dcal_X$ are associative, compatible to the units, and satisfies the exchange law but this is completely trivial when we write the corresponding co-limits expressing these compositions, and the relation giving the value of $\pi^{\epsilon}_k(a \#_n b)$ follow also directly from the definitions.
}
}

\block{\label{Const_MorphDxDstar}One can construct two morphisms $X^* \rightarrow \Dcal_X$ and $\Dcal_X \rightarrow X^*$ as follow:

To any $(v,\alpha,\lambda) \in \Dcal_X$ one can associate $\lambda^*(\alpha) \in X^*$. It is immediate that this is a morphism of $\infty$-categories.

$X^* \rightarrow \Dcal_X$ is constructed as follow: on the underlying $n$-category it is the natural isomorphism between their underlying $n$-categories mentioned in proposition \ref{DefOfDcalX}. Because of the universal property of $X^*$ we only need to say what are the images of the $(n+1)$-generators of $X^*$. For each $(n+1)$-generator $x$ of $X^*$ there is a $J_{n+1}$-plex $\underline{a}$ of dimension $n+1$ and a map $\chi_x: \underline{a} \rightarrow X$ corresponding to $x$. One defines the image of $x$ in $\Dcal_X$ as the triple $(\underline{a},a, \chi_x)$, which can easily be checked to have the correct source and target.

Finally, it is immediate to check that the composite $X^* \rightarrow \Dcal_X \rightarrow X^*$ is the identity of $X^*$, and that these two constructions are functorial in $X$. One hence has in particular proved that:

\Prop{$X^*$ is a retract of $\Dcal_X$, functorially in $X$.}

}

\block{\label{Def_Jnp1Pasting}\Def{A $J_{n+1}$-polyplex is a pair $(v,\alpha)$ with $v \in J_{n+1}$ and $\alpha \in v^*$ whose isomorphism class is in the image of the morphism $(J_{n+1}1)^* \rightarrow \Dcal$ constructed in \ref{Const_MorphDxDstar}.

}

We will also denote such pairs $(\underline{x},x)$ following our conventions, even if these are not exactly polyplexes in the usual sense.

Note that as $(J_{n+1} 1)^*$ is a retract of $\Dcal$, this already shows that isomorphism classes of $J_{n+1}$-polyplexes are exactly indexed by arrows of $(J_{n+1} 1)^*$ and form a sub-$\infty$-category of $\Dcal$. However this time we do not want to systematically identify polyplexes with arrows of $(J_{n+1} 1 )^*$ as the polyplex attached to an arrow of $J_{n+1}$ is only well defined up to non-unique isomorphism. We might still occasionally use this identification as a notational shortcut.}

\block{\label{Dicuss:condCforJn+1Polyplexes}One also have, essentially by construction, a version of condition $(C)$ for these newly defined $J_{n+1}$-polyplexes: given two composable arrows of $(J_{n+1} 1)^*$ corresponding to $J_{n+1}$-polyplexes $(\underline{f},f)$ and $(\underline{g},g)$, then because the map $(J_{n+1} 1)^* \rightarrow \Dcal$ has been constructed as a morphism of $\infty$-categories, the $J_{n+1}$-polyplex corresponding to their composition, is the composite of $(\underline{f},f)$ and $(\underline{g},g)$ in $\Dcal$, i.e. it is indeed $(\underline{f} \coprod_{\underline{p}} \underline{g} ,f \#_k g)$ as condition $(C)$ requires it.

}

\block{\label{Prop_ImageX*inDx}\Prop{The image of $X^*$ in $\Dcal_X$ is exactly the set of (equivalence classes of) triples $(v,\alpha,\lambda)$ such that $(v,\alpha)$ is a $J_{n+1}$-polyplex.}

\Dem{First we observe that the set $W$ of $(v,\alpha,\lambda)$ such that $(v,\alpha)$ is a polyplex is indeed an $\infty$-category. Indeed as the map $(J_{n+1}1)^* \rightarrow \Dcal$ is a monomorphism (it has a retraction) its image is a sub-$\infty$-category, and hence its pre-image along the morphism $\Dcal_X \rightarrow \Dcal$ is also an $\infty$-category, and this is exactly the subset $W$.

This implies that $X^*$ is included in $W$, indeed, all the arrows of dimension $ \leqslant n$ and all the generators of dimension $n+1$ are sent into $W$, so one has indeed $X^* \subset W$.

Finally we will prove by induction on the arrows of $(J_{n+1}1)^*$ that if $x$ is an arrow of $(J_{n+1}1)^*$ corresponding to a polyplex $(\underline{v},v)$ and if $\lambda$ is any morphism from $\underline{v}$ to $X$ then $(\underline{v},v,\lambda) \in X^* \subset \Dcal_X$.

If $\underline{x}$ is a generator of $(J_{n+1}1)^*$, then the corresponding polyplex is $(\underline{x},x)$, a morphism $\lambda: \underline{x} \rightarrow X$ just specify a cell of $X$ and $(\underline{x},x,\lambda)$ is the image in $\Dcal_X$ of the corresponding cell of $X$.

Assume now that $x =a \#_k b$ in $(J_{n+1} 1)^*$ and $a$ and $b$ satisfies the induction hypothesis. $a$ and $b$ corresponds respectively to polyplexes $(\underline{f},f)$ and $(\underline{g},g)$ and $x$ corresponds to $(\underline{f} \coprod_{\underline{p}} \underline{g}, f \#_k g )$ where $\underline{p}$ is the $J_n$-polyplex corresponding to the $k$-source of $a$ and the $k$-target of $b$. If $\lambda : \underline{f} \coprod_{\underline{p}} \underline{g} \rightarrow X$ is any map then it is exactly the data of two maps $\lambda_1 : \underline{f} \rightarrow X$ and $\lambda_2 : \underline{g} \rightarrow X$ which agrees on $\underline{p}$. Both $(\underline{f},f,\lambda_1)$ and $(\underline{g},g,\lambda_2)$ are arrows of $\Dcal_X$ which are in the image of $X^*$ by the induction hypothesis, and $x$ is exactly their $k$-composite, so it also belongs to the image of $X^*$.
}
}

\block{\label{Th_mainGoodClass}\Th{Let $J$ be a class of polygraphs such that $J_n$ is a good class of polygraphs satisfying condition $(C)$ of \ref{Def_CondC}. Then for $X \in J_{n+1}$ one has:

\[ X^* = \coprod_{\underline{v} \text{ a }J_{n+1}\text{-polyplex}} \left( Hom(\underline{v},X)/G_{\underline{v}} \right) \]

If all the $G_{\underline{v}}$ are trivial then $J_{n+1}$ is a good class of polygraphs satisfying condition $(C)$, with the two notions of $J_{n+1}$-polyplexes (from \ref{Prop_goodclass} and \ref{Def_Jnp1Pasting}) being equivalent. Conversely, if $J_{n+1}$ is a good class of polygraphs then all the $G_{\underline{v}}$ are trivial.
}

\Dem{We have seen in \ref{Prop_ImageX*inDx} that $X^*$ can be described as the subset of $\Dcal_X$ of equivalence classes of triples $(\underline{v},v,\lambda)$ where $(\underline{v},v)$ is a polyplex. Two such triples are equivalent exactly if the polyplexes are isomorphic (as polyplex) and if $\lambda$ and $\lambda'$ are conjugate under the action of $G_{\underline{v}}$, and this already proves that $X^*$ has the given description. If all the $G_v$ are trivial then this proves that $X \mapsto X^*$ is familially representable. We have seen in \ref{Prop_GoodnImpPrenpp} that $J_{n+1}$ is an effective presheaf category hence $J_{n+1}$ is a good class of polygraphs. The explicit formula for composition in $\Dcal_X$ given in \ref{DefOfDcalX} shows that composition of two polyplexes is given by a pushout (see \ref{Dicuss:condCforJn+1Polyplexes}), and hence that $J_{n+1}$ satisfies condition $(C)$. Conversely, if $J_{n+1}$ is a good class of polygraphs, it means that the functor $X \mapsto X^*$ is familially representable:

\[ X^* \simeq \coprod_{i \in I} Hom(w_i,X) \] 

By comparing this formula to the one we gave in the statement of the theorem and looking at the value of this functor on the terminal $J_{n+1}$-polygraph and on the map to the terminal $J_{n+1}$-polygraph one concludes that there is a bijection between the $w_i$ and the $J_{n+1}$-polyplexes as they were defined in this section such that $Hom(w_i,X) \simeq Hom(\underline{v},X)/G_{\underline{v}}$. Lemma \ref{Lem_QuotientRep} below concludes the proof.

}

}

\block{\label{Lem_QuotientRep}\Lem{Let $\Ccal$ be any category, $v$ an object of $\Ccal$ endowed with an action of an ordinary group $G$. Then the set valued functor:

\[ X \mapsto Hom(v,X)/G \]

is representable if and only if the action of $G$ on $v$ is the trivial action.

}

\Dem{If $G$ acts trivially, then the functor is representable by $v$. Conversely, assume that $w$ is an object such that one has an isomorphism $\theta_X : Hom(v,X)/G  \overset{\simeq}{\rightarrow} Hom(w,X)$ functorial in $X$. Let $i:w \rightarrow v$ be the element $\theta_v(1_v)$ and let $p:v \rightarrow w$ be any element of the $G$-class $\theta_w^{-1}(1_w) \in Hom(v,w)/G$.

By applying to $1_v \in Hom(v,v)/G$ the naturality square:

\[\begin{tikzcd}[ampersand replacement=\&]
Hom(v,v)/G \arrow{r}{p \circ } \arrow{d}{\theta_v} \&  Hom(v,w)/G \arrow{d}{\theta_w} \\
Hom(w,v) \arrow{r}{p \circ} \& Hom(w,w)
\end{tikzcd}
\]

one gets : $ p \circ i = 1_w $

By applying to $1_w \in Hom(w,w)$ the naturality square:

\[\begin{tikzcd}[ampersand replacement=\&]
Hom(w,w) \arrow{r}{i \circ } \arrow{d}{\theta_w^{-1}} \& Hom(w,v) \arrow{d}{\theta_v^{-1}} \\
Hom(v,w)/G \arrow{r}{i \circ } \&  Hom(v,v)/G
\end{tikzcd}
\]

One obtains that $i \circ p$ is equal to $1_v$ in $Hom(v,v)/G$, i.e. it is in the image of $G$, and hence it is invertible. This implies that $i$ and $p$ are inverse of each other.

Finally, by applying to $1_v \in Hom(v,v)/G$ the following naturality square for any $g \in G$:

\[\begin{tikzcd}[ampersand replacement=\&]
Hom(v,v)/G \arrow{r}{g \circ } \arrow{d}{\theta_v} \&  Hom(v,v)/G \arrow{d}{\theta_v} \\
Hom(w,v) \arrow{r}{g \circ} \& Hom(w,v)
\end{tikzcd}
\]

one gets that $g \circ i =i$, and as $i$ is an isomorphism this implies that $g=Id$ as an automorphism of $v$, i.e. that the action of $G$ on $v$ is trivial.

}
}

\block{\label{Cor_CondCforGoodClass}\Cor{Condition $(C)$ of \ref{Def_CondC} is satisfied by any good class of polygraphs.}

\Dem{Let $J$ be a good class of polygraphs, we prove by induction on $n$ that $J_n$ satisfies condition $(C)$. Condition $(C)$ is clearly satisfied for the class of all $1$-polygraphs, so $J_1$ satisfies it. If $J_n$ is a good class of polygraphs satisfying condition $(C)$ then one can apply theorem \ref{Th_mainGoodClass} to it. As $J_{n+1}$ is a good class of polygraphs all the $G_{\underline{v}}$ are trivial and hence $J_{n+1}$ is a good class of polygraphs satisfying condition $(C)$. As condition $(C)$ only involves finite polygraphs in $J$ it is enough to check it for $J_n$ for all $n$, hence this concludes the proof.}
}

\subsection{Automorphisms of polyplexes and the case of source-positive polygraphs}
\label{subsection_AutomPD-SPpolyg}

\blockn{In this section we will study in more details the automorphisms of a polyplex in the sense of definition \ref{Def_AutPD}, i.e. automorphisms of the underlying polygraph whose action on the free $\infty$-category it generates preserves the universal cell. We are still working under the same assumption as in the previous section: $J$ is a class of polygraphs such that $J_n$ is a good class of polygraphs, and we are interested in the $J_{n+1}$-polyplexes as they were defined in \ref{Def_Jnp1Pasting}, and more precisely in their automorphisms. Our goal is to find criterions under which these automorphisms are all trivial. The main result being theorem \ref{Th_SourcePositiveIsGood}, showing that in the case where $J$ is the class of source-positive polygraphs then all these automorphisms are indeed trivial and hence the class of all source-positive polygraphs is a good class of polygraphs. 

All the polyplexes of dimension $\leqslant n$ are in the good class $J_n$ and hence have no automorphisms. Hence we are only interested in the polyplexes of dimension $n+1$. 
}

\block{\label{disc:Cond(C)JN+1}By \ref{Cor_CondCforGoodClass}, $J_n$ satisfies ``condition $(C)$'' of \ref{Def_CondC}. We also recall that it follows from the results of the previous section (see \ref{Dicuss:condCforJn+1Polyplexes}), that despite not being known to be a good class of polygraphs yet, the class $J_{n+1}$ also satisfies a form of condition $(C)$, in the sense that if two arrows $u,v \in X^*$ for $X \in J_{n+1}$ are represented by polyplexes $\chi_u:\underline{p} \rightarrow X$ and $\chi_v:\underline{q} \rightarrow X$ then their composite $u \#_k v$ (if it exists) will be represented by a polyplex isomorphic to $\underline{p} \coprod_{\underline{c}} \underline{q}$ where $\underline{c}$ is the polyplex representing the arrow $\pi^+_k u =\pi^-_k v$. This is because the isomorphisms of theorem \ref{Th_mainGoodClass} showing that any arrow in a $J_{n+1}$-polygraph $X^*$ can be represented by a unique polyplex $\underline{p} \rightarrow X$ up to non-unique isomorphism, is not just an isomorphism of sets, but an isomorphism of $\infty$-categories between $X^*$ and the subobject of the $\infty$-category $\Dcal_X$ (constructed in \ref{DefOfDcalX}) corresponding to polyplexes (see proposition \ref{Prop_ImageX*inDx}), hence one can compute the polyplexes representing a composition in $X^*$ in terms of the composition in $\Dcal_X$, which is given by taking a pushout.

It is also the case that the polyplex representing a generator is just the corresponding plex, which is constructed (see \ref{Prop_GoodnImpPrenpp}) exactly as in a good classes of polygraphs, by gluing the polyplexes corresponding to its source and target along their common boundary and adding a unique cell of maximal dimension. This is because this is how the morphism $X^* \rightarrow \Dcal_X$ is constructed in the first place.
}

\block{If $X$ is a $J_{n+1}$-polygraph and $v \in X^*$ is an arrow of dimension $\leqslant n+1$ then because of theorem \ref{Th_mainGoodClass} there is a polyplex $(\underline{p},p)$ and a polygraphic morphism $\chi_v : \underline{p} \rightarrow X$ such that $\chi_v^*(p)=v$, moreover $\chi_v$ and $\underline{p}$ are unique up to non-unique isomorphism. In particular the image of the morphism $\chi_v$ in $X$ is well defined. We call this subpolygraph of $X$ the support of $v$ and we denote it $Supp(v)$. One can check that this agree with Makkai's definition of the support of an arrow given in lemma $5.(5)$ of \cite{makkai2005word}, indeed:

\Lem{The support of an arrow can be defined inductively by the formula:

\begin{itemize}

\item $Supp(x) = \{x \}$ if $x$ is a $0$-cell of $X$.

\item $Supp(x)=\{x\} \cup Supp( \pi^-_{n-1} x) \cup Supp( \pi^+_{n-1} x)$ if $x$ is a $n$-cell of $X$.

\item $Supp(a \#_k b) = Supp(a) \cup Supp(b)$.

\end{itemize}

Moreover $Supp(v)$ is the smallest subpolygraph $Y \subset X$ such that $v \in Y^*$

}

\Dem{The induction formula follows immediately from the description of the polyplex representing a generator or a composite given in the discussion in \ref{disc:Cond(C)JN+1}. For the last remark $Supp(v)$ is a polygraph because it is the image of a polygraph, and if $v \in Y^*$ for $Y \subset X$ then the polyplex representing $v$ has to factor into $Y$ and hence $Supp(v) \subset Y$.
}
}

\block{In this section we will say that a cell $x \in X$ appears in an arrow $v\in X^*$ if $x \in Supp(v)$. Because of the previous lemma $x$ appears in a composite $f \#_k g$ if and only if $x$ appears in $f$ or in $g$, and if $x$ appears in $f$; and if $f$ belong to $Y^*$ for some subpolygraph $Y \subset X$ then $x \in Y$. For an $n$-cell $x$ and an $n$-arrow $f$ one has $x \in Supp(f)$ if and only if the ``number of time $x$ appears in $f$'' (in the sense of the third point of \ref{Discus:SyntacticalPropOfPolygraphs}) is non-zero, so the two uses of ``$x$ appears in $f$'' are compatible.
}

\block{\label{Lem:injectivityOnKandK-1cells} Even in a good class, polyplexes can have surprisingly pathological behavior, typically the gluing used for composition, or the boundary inclusion can fail to be inclusion, see for example \ref{ExIsomorphicPolyplex}. But it appears that, as the following lemma shows, these problems only appears in codimension at least $2$. Most of our arguments in this section will be at the level of $n$-cells and $n+1$-cells in $n+1$-polyplexes, so will avoid these problems.
 
\Lem{
\begin{itemize}

\item Let $\underline{a}$ be a $(n+1)$-polyplex, then the two maps $\pi^{\epsilon}_{n} \underline{a} \rightarrow \underline{a}$ are injective on $n$-cells.

\item Let $\underline{a}$ and $\underline{b}$ be a pair of $k$-composable polyplexes of dimension at most $(n+1)$, then the two maps $\underline{a},\underline{b} \rightrightarrows \underline{a} \#_k \underline{b}$ are injective on $n$-cells and $(n+1)$-cells.

\end{itemize}
}

\Dem{We prove the first point by induction on the polyplex $\underline{a}$, more precisely, we prove by induction on arrows of $(J1)^*$ that the corresponding polyplex $\underline{a}$ satisfies the first point. The second point will appear as a byproduct of this proof.

If $\underline{a}$ is a $(n+1)$-plex, then it is constructed by adding a single generator to the gluing of its source and target along their $n-1$-dimensional boundary. In particular the map:

\[ \pi^-_{n} \underline{a} \coprod_{\pi^-_{n-1} \underline{a} \coprod \pi^+_{n-1} \underline{a} } \pi^+_{n} \underline{a} \hookrightarrow \underline{a} \]

is always a monomorphism. Now the functor sending a polygraph to its set of $n$-cells preserves colimits, and the polygraph $\pi^-_{n-1} \underline{a} \coprod \pi^+_{n-1} \underline{a} $ has no $n$-cells as it is of dimension $n-1$, so the set of $n$-cells of the left hand side above is just the disjoint union of the $n$-cells of $\pi^-_{n} \underline{a}$ and $\pi^+_{n} \underline{a}$, which proves the claim for a plex.

We now consider a $(n+1)$-polyplex $\underline{v} = \underline{a} \#_k \underline{b}$ with $\underline{a}$ and $\underline{b}$ satisfying the induction hypothesis. If $k \geqslant n+1$ then $\underline{v} = \underline{a} = \underline{b}$ and everything is trivial, so one can freely assume that $k \leqslant n$.

We first prove that this specific composition satisfies the second point of the lemma. We need to treat two cases separately: if $k<n$ then $\pi^+_k \underline{a}$ as no $n$-cells and no $(n+1)$-cells, so the set of $n$ and $(n+1)$-cells of $\underline{v} = \underline{a} \coprod_{\pi^+_k \underline{a}} \underline{b}$ is just the disjoint union of those of $\underline{a}$ and $\underline{b}$, so the maps $\underline{a}, \underline{b} \rightrightarrows \underline{v}$ are injective on $n$ and $(n+1)$-cells. If $k=n$ the same argument show injectivity on $(n+1)$-cells, our induction hypothesis also claim that the two maps $\pi^+_{n} \underline{a} \rightrightarrows \underline{a},\underline{b}$ are injective on $n$-cells. The set of $n$-cells of $\underline{v}$ is hence obtained as a pushout of two monomorphisms, which proves that as claimed, the two structural maps $\underline{a},\underline{b} \rightrightarrows \underline{v}$ of this pushout are injective on $n$-cells.

We now finish the proof of the first point for $\underline{v}$. We also need to distinguish two cases. If $k=n$ then $\pi^-_{n} \underline{v} = \pi^-_{n} \underline{a}$ hence one can consider the composite:

\[ \pi^-_{n} \underline{v} = \pi^-_{n} \underline{a} \rightarrow \underline{a} \rightarrow \underline{v} \]

The first map is injective on $n$-cells by the induction hypothesis and the second map is also injective on $n$-cells because of the second point of the lemma (which we just proved in this case). The case of target is obtained in the exact same way using $\pi^+_n \underline{v} = \pi^+_n \underline{b}$ instead.

If $k< n$, then $\pi^{\epsilon}_{n}(\underline{a} \#_k \underline{b}) = \pi^{\epsilon}_{n}\underline{a} \#_k \pi^{\epsilon}_{n} \underline{b}$, one considers the diagram:

\[\begin{tikzcd}[ampersand replacement=\&]
\pi^{\epsilon}_{n} \underline{a} \arrow{d} \& \arrow{l} \pi^+_k \underline{a} \arrow{d} \arrow{r} \& \pi^{\epsilon}_{n} \underline{b} \arrow{d} \\
\underline{a} \& \arrow{l} \pi_k^+ \underline{a} \arrow{r} \& \underline{b} \\
\end{tikzcd}\]

the map $\pi^{\epsilon}_{n} \underline{v} \rightarrow \underline{v}$ is the comparison map between the pushout of the top line and the pushout of the bottom line. Now, the middle object have no $n$-cells and the two other vertical maps are injective on $n$-cells, so this comparison map is injective on $n$-cells.
}
}

\block{\label{Lem_LocCellInPD}
\Lem{Let $(\underline{p},p)$ be a $J_{n+1}$-polyplex of dimension $n+1$ and $e$ a cell of $\underline{p}$ of dimension $\leqslant n$ then one has either:

\begin{itemize} 

\item $e$ appears in $\pi^-_n(p)$.

\item $e$ appears in $\pi^+_n(x)$ for $x$ an $(n+1)$-cell of $\underline{p}$.

\end{itemize}

Moreover, if $e$ is of dimension $n$ exactly then the two possibilities above cannot occur simultaneously and in the second case the $(n+1)$-generator $x$ such that $e$ appears in $\pi^+_n(x)$ is unique.

}

One also has a dual statement obtained by exchanging source and target everywhere whose proof is exactly the same.

\Dem{We will prove the lemma by induction on arrows of $(J_{n+1}1)^*$, corresponding to the polyplex $\underline{p}$. 

For a generator of $(J_{n+1}1)^*$, the polyplex $(\underline{p},p)$ is a $J_{n+1}$-plex with its unique $(n+1)$-cell. As mentioned in proposition \ref{Prop_GoodnImpPrenpp}, a $J_{n+1}$-plex of dimension $n+1$ is always constructed by gluing two parallel $J_n$-polyplexes $\pi^-_n \underline{p}$ and $\pi^+_n \underline{p}$ along their common boundary and adding the $(n+1)$-cell $p$ with the $\pi^{\epsilon}_n p $ being the images of universal cells of the $\pi^{\epsilon}_n \underline{p}$.  Hence any generator $e$ of dimension $\leqslant n$ is either in the image of $\pi^-_n \underline{p}$, which corresponds to the first case of the dichotomy in the lemma, or in the image of $\pi^+_n \underline{p}$ which corresponds to the second case of dichotomy. Moreover if $e$ is of dimension $n$, then as $\pi^-_n \underline{p}$ and $\pi^+_n \underline{p}$ are glued together along a polygraph of dimension $<n$, $e$ can only belong to the image of one of them and so the two possibilities are indeed incompatible.

For a composite $\underline{a} \#_k \underline{b}$ in $(J_{n+1}1)^*$, it corresponds to a polyplex of the form:

\[ (\underline{a} \coprod_{\underline{c}} \underline{b}, a \#_k b) \]

where $\underline{c} = \pi^+_k \underline{a} = \pi^-_k \underline{b}$, in particular it is of dimension $\leqslant k$.

We need to distinguish several cases, depending on the dimension of $\underline{a}$ and $\underline{b}$ and the value of $k$:

We remind that because of lemma \ref{Lem:injectivityOnKandK-1cells} the maps $\underline{a},\underline{b} \rightrightarrows \underline{a} \#_k \underline{b}$ are injective on $n$-cells and $n+1$-cells, this is a key point in the inductive proof of the ``uniqueness part'' of the lemma, which only involves $n$ and $n+1$-cells. Because of this fact, we will tend to not distinguish between the $n$ and $(n+1)$-cells of $\underline{a}$ and $\underline{b}$ and their images in the composite.

\begin{itemize}

\item If $k=n$ and both $\underline{a}$ and $\underline{b}$ are of dimension $n+1$, then the cell $e$ appears in either $a$ or $b$, i.e. is either in the image of $\underline{a}$ or in the image of $\underline{b}$. If it is in the image of $\underline{a}$, then by induction (a pre-image of) $e$ appears either in the source of $a$, in which case it appears also in the source of $a \#_n b$, or it appears in the the target of a $(n+1)$-cell $x$ of $\underline{a}$, which is also a $(n+1)$-cell of $\underline{a}\#_k \underline{b}$. If now $e$ is in the image of $\underline{b}$, then by induction either $e$ appears in the source of $b$, which coincide with the target of $a$, hence $e$ also belong to the image of $\underline{a}$ and we are brought back to the previous case. Or $e$ appears in the target of a $(n+1)$-cell $x$ of $\underline{b}$, which is also a $(n+1)$-cell of the composite.

We now assume that the cell $e$ is of dimension $n$ exactly and we prove the uniqueness part of the result.

If $e$ appears in both $a$ and $b$, it means that it is both in the image of $\underline{a}$ and $\underline{b}$ in:

 \[ \underline{a} \#_k \underline{b} = \underline{a} \coprod_{\pi^+_n \underline{a}} \underline{b} \]

which means that it has to belong to the image of $\pi^+_n \underline{a} = \pi^+_n \underline{b}$, but if it appears in the source of $b$ it cannot appears in the target of any $(n+1)$-generator of $b$ by the induction hypothesis, so by the induction hypothesis for $a$ it (exclusively) either appears in the source of $a$ or in the target of a unique $(n+1)$-generator in $a$ and this proves the uniqueness part in this case. If it appears in only one of $a$ and $b$, then one gets the uniqueness directly from the inductive application of the uniqueness result in $a$ and $b$ separately (it cannot appear in the source of $b$ in this case).

\item If both $a$ and $b$ are of dimension $n+1$ but $k<n$: a cell $e$ appears in either $a$ or $b$. In both case it appears either the target of an $(n+1)$-cell or the source of $a$ or $b$, but in this case both the support of the source of $a$ and the support of the source of $b$ are included in the support of the source of $a \#_k b$ so it concludes the proof immediately. If $e$ is of dimension $n$ exactly then as $\underline{a}$ and $\underline{b}$ are glued together on something of dimension $<n$, the generator $e$ cannot be in both of them simultaneously and so one gets the uniqueness result by simply applying the uniqueness result inductively in $a$ and $b$ separately.

\item If $k<n$ and either $a$ or $b$ is of dimension $\leqslant n$, for example assume $b$ is of dimension $\leqslant n$ (the proof is exactly the same if $a$ is of dimension $\leqslant n$). In this situation $\pi^-_n(a \#_k b) = \pi^-_n(a) \#_k b$ hence any cell $e$ that appears in $b$ also appears in $\pi^-_n(a \#_k b)$ which proves the claim for such cells. For any cell $e$ appearing in $a$, one can just apply the induction hypothesis inside $\underline{a}$: either (a preimage of) $e$ appears in the source of $a$ or in the target of some $(n+1)$-cell of $a$, in both case it proves the claim. If $e$ is of dimension $n$ exactly then $e$ cannot be both in the image of $\underline{a}$ and $\underline{b}$ simultaneously as they are glued along a polygraph of dimension $k<n$, if $e$ appears in $b$ then it appears in the source of $a \#_k b$ and there is no $(n+1)$-cell it could appears in because all the $(n+1)$-cells are in $a$, and if it appears in $a$ then it appears (uniquely) either in the source of $a$ or in the target of a unique $(n+1)$-generator of $a$ and this concludes the proof of this case.

\item All the other cases are trivial: if $k=n$ and $a$ and $b$ are not both of dimension $n+1$ then the composite is equal to $a$ or $b$, if neither $a$ nor $b$ is of dimension $n+1$ then the composite is not of dimension $n+1$ and this also concludes the proof.

\end{itemize}

}
}

\blockn{The next lemma is our key tool to show that the automorphisms of a given polyplex are trivial. It shows that any automorphism have some fixed points, and that being a fix point of an automorphisms tend to be ``contagious''.}

\block{\label{Lem_wFixContagious}\Lem{Let $(\underline{p},p)$ be a $J_{n+1}$-polyplex of dimension $n+1$ and $w \in G_{\underline{p}}$ be an automorphism. 

\begin{enumerate}

\item If an arrow $t \in \underline{p}^*$ of dimension $\leqslant n$ is fixed by $w^*$, then all the cells that appears in $t$ are fixed by $w$.

\item $w$ fixes all the cells that appears in $\pi^-_n(p)$ or $\pi^+_n(p)$.

\item If $w$ fixes some $(n+1)$-cell $x$ then it also fixes all the cells that appears in $\pi^-_n(x)$ and $\pi^+_n(x)$.

\item If $w$ fixes an $n$-cell $a$ that appears in the source or the target of an $(n+1)$-cell $x$, then $x$ is also fixed by $w$.

\end{enumerate}

 }

\Dem{
\begin{enumerate}

\item let $t \in \underline{p}^*$ be an arrow fixed by $w^*$. As $t$ is of dimension smaller than $n$ it belongs to the $\infty$-category generated by $n$-th skeleton of $\underline{p}$ which is a $J_n$-polygraph. In particular $t$ is uniquely represented by a $J_n$-polyplex $\chi_t:\underline{v} \rightarrow \underline{p}$, in the sense that $\chi_t^*(v)=t$. As $w^* (t) =t$ one has that $w \circ \chi_t $ also represents $t$, so by uniqueness of the representation $w \circ \chi_t =  \chi_t$ and hence all cells that are in the image of $\chi_t$ are fixed by $w$. Those are exactly the cells that appears in $t$.

\item As $w \in G_{\underline{p}}$, one has $w^*(p)=p$, and so $w$ also fixes $\pi^-_n(p)$ and $\pi^+_{n}(p)$. As those are arrows of dimension $\leqslant n$, one can apply point $1.$ to conclude.

\item If $w$ fixes an arrow $x$ then it also fixes the $\pi^{\epsilon}_n(x)$, which are arrows of dimension $\leqslant n$, hence the point $1.$ applies and $w$ fixes all the cells that appears in these.

\item If a cell $a$ appears in $\pi^{\epsilon}_n(x)$ for some $\epsilon$ and some arrow $x$, then $w(a)$ will appears in $\pi^{\epsilon}_n(w(x))$. Hence if $w(a)=a$ it means that $a$ appears in both $\pi^{\epsilon}_n(x)$ and $\pi^{\epsilon}_n(w(x))$, but as $a$ is of dimension $n$ the uniqueness part of lemma \ref{Lem_LocCellInPD} (or its dual in the case $\epsilon=-$) shows that $x=w(x)$.

\end{enumerate}

}
}

\block{\label{Prop_AutActsOnTopCellsOnly}\Prop{If an automorphism $w$ of a polyplex fixes all the $(n+1)$-cells then it is the identity.}

In particular, $G_{\underline{v}}$ can be identified with a sub-group of the group of permutations of the $(n+1)$-cells of $\underline{v}$.

\Dem{It follows immediately from lemma \ref{Lem_LocCellInPD} together with points $2.$ and $3.$ of lemma \ref{Lem_wFixContagious}.
}

}

\block{\label{Th_SourcePositiveIsGood}\Th{The class of all source-positive polygraphs is a good class of polygraphs.}

\Dem{
Let $J$ be the class of all source-positive polygraphs. We proceed by induction using theorem \ref{Th_mainGoodClass}: $J_1$ is a good class of polygraphs (it is the category of graphs). We assume by induction that $J_n$ is a good class of polygraphs, and we need to show that for any $J_{n+1}$-polyplex $(\underline{p},p)$ of dimension $n+1$, any automorphism $w \in G_{\underline{p}}$ is the identity.

\bigskip

The general idea is to use lemma \ref{Lem_wFixContagious} to propagate the fact that $w$ is the identity from the source of $\underline{p}$ to all the $(n+1)$-cells through the $n$-cells connecting them, and then use lemma \ref{Prop_AutActsOnTopCellsOnly} to conclude that $w$ is the identity.

In order to make this formal, one introduces the following notion: A chain of cells in $\underline{p}$ is a finite sequence $x_1,\dots,x_k$ of $(n+1)$-cells of $\underline{p}$ such that:

\begin{itemize}

\item There is an $n$-cell which appears both in the source of $x_1$ and in $\pi^-_n(p)$.

\item For each $i>0$, there is an $n$-cell which appears both in the source of $x_{i+1}$ and in the target of $x_i$.

\end{itemize}

Lemma \ref{Lem_wFixContagious} shows that any such chain is fixed by $w$: indeed the $n$-cell appearing in the source of $p$ is fixed because of point $2.$ of \ref{Lem_wFixContagious} and this implies that $x_1$ is fixed because of point $4.$, and then inductively if $x_i$ is fixed then the $n$-cell which appears in both the target of $x_i$ and the source of $x_{i+1}$ is fixed because of point of $3.$ and this implies that $x_{i+1}$ is fixed by point $4.$ of the lemma.

\bigskip

We will now show that any $(n+1)$-cell of $\underline{p}$ belongs to such a chain. More precisely, we will prove by induction on arrows of $(J_{n+1}1)^*$ that in the corresponding polyplex $(\underline{p},p)$ any $(n+1)$-cell of $\underline{p}$ appears in such a chain.

If $(\underline{p},p)$ corresponds to a generator of $(J_{n+1} 1 )^*$, i.e. if $\underline{p}$ is a $J_{n+1}$-plex, then the unique $(n+1)$-dimensional cell is the $x_1$ of a chain exactly because one has assumed that its source contains at least one $n$-cell.

If $(\underline{p},p)$ corresponds to a composite $\underline{a} \#_k \underline{b}$ in $(J_{n+1}1)^*$, it is of the form:

\[ (\underline{a} \coprod_{\underline{c}} \underline{b}, a \#_k b) \]

Let $x$ be an $(n+1)$-cell in this polygraph $\underline{p}$. Then $x$ belongs to either the image of $\underline{a}$ or the image of $\underline{b}$. As the claim is only about $n$-cells and $(n+1)$-cells and lemma \ref{Lem:injectivityOnKandK-1cells} shows that the maps $\underline{a},\underline{b} \rightrightarrows \underline{a} \#_k \underline{b} = \underline{p}$ are injective on $n$ and $(n+1)$-cells one can simply identifies $\underline{a}$ and $\underline{b}$ with their image in $\underline{p}$ and just says that $x$ belongs to $\underline{a}$ or $\underline{b}$.

 If $x$ is in $\underline{a}$ it belongs to a chain in $\underline{a}$ which starts on an $n$-cell of the source of $\underline{a}$, and the source of $\underline{a}$ is always included in the source of $\underline{a} \#_k \underline{b}$. If $x$ is in $\underline{b}$, there is a chain $(x_1,\dots,x)$ in $\underline{b}$, in particular there is an $n$-cell $e$ which appears both in the source of $x_1$ and in the source of $b$. If $k<n$, then the source of $\underline{b}$ is included in the source of $\underline{p}$ and $(x_1,\dots,x)$ is already the chain we are looking for. If $k=n$ then the source of $b$ is exactly the target of $a$ (in $\underline{p}$) and by lemma \ref{Lem_LocCellInPD} $e$ appears either in the source of $a$ (in which case the proof is finished) or in the target of some $(n+1)$-cell $y$ of $\underline{a}$ in which case there is a chain $y_1,\dots,y_i=y$ in $\underline{a}$ and $y_1,\dots,y_n,x_1,\dots,x$ is a chain in $\underline{p}$ which contains $x$. Indeed, the only ``chain condition'' that does not follow from the fact that $(x_i)$ and $(y_i)$ are respectively chains in $\underline{b}$ and $\underline{a}$ is the one between $y_n=y$ and $x_1$ and for this one, the element $e$ provides the link as it appears both in the source of $x_1$ and in the target of $y$.

}

}

\block{\Cor{The class of positive polygraphs and the class of opetopic polygraphs are both good classes of polygraphs. In particular they are presheaf categories.}

\Dem{They are both included in the class of target-positive polygraphs, which is a good class of polygraphs exactly as the class of source-positive polygraphs. Hence by proposition \ref{Prop_goodclass} they are also good classes of polygraphs.}

}

\subsection{Final remarks and (counter)examples}

\block{In this subsection one just make some general comment and give some example of surprising behavior one can have, even in good classes of polygraphs. We also mention some further questions that would be interesting:

\begin{itemize}

\item Can we use the inductive description of plexes to show that the ``opetopic plexes'' are the same as the opetopes ? This would give a direct proof that opetopic polygraphs are the same as opetopic sets. Note: shortly after the publication of this paper, Cedric Ho tanh has given in \cite{thanh2018equivalence} a direct proof of the equivalence between opetopic polygraphs and opetopic sets which does not relies on our approach.

\item Can we give a more convenient description of the positive plexes or the source-positive plexes ? We are after some generalization of the various notion of $n$-pasting diagrams that have been devised by M.Johnson \cite{johnson1989combinatorics}, A.J.Power \cite{power1991n}, R.Steiner \cite{steiner2004omega} or R.Street \cite{street1991parity}. This has been achieved in \cite{hadzihasanovic2018combinatorial} by Amar Hadzihasanovic for a large class of polygraphs that he calls the regular\footnote{They are only a subclass of what I called regular polygraphs in \cite{henry2018regular}, but the two notions are very close and we believe a slight modification of his description can encompass the slightly more general notion of regular polygraphs of \cite{henry2018regular}.} polygraphs.

\end{itemize}

}

\block{\label{ExIsomorphicPolyplex}We now give a promised examples of two positive $3$-polyplexes with isomorphic underlying polygraphs, but different universal arrows. It also contains some examples of $3$-plexes whose ``boundary maps'' are not injective. We also believe that this example shows that the second question above is hard: the composition formed by these two polyplexes uses the exact same variables (meaning they have the same underlying polygraphs) but produces different results which hence only differs by the 'order' in which certain cells are composed. But the key point in all the theory of pasting diagrams mentioned is to restrict to situations where there is a unique possible composition order or where the composition order does not matter, and the following example seems to be far out of the scope of such situations. This is not a clear counter example to the existence of a nice description, but clearly an example we should meditate before attacking this question. 

We start with the following three $2$-polyplexes:

 \[ \begin{array}{c | c | c} 
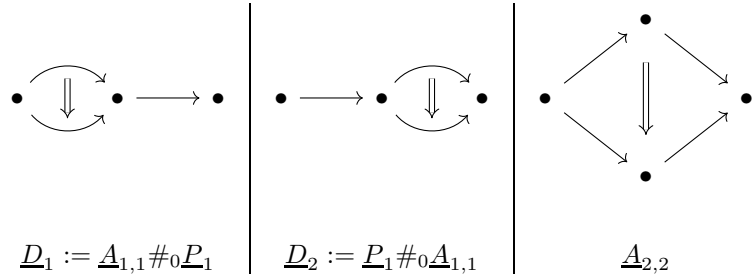
\begin{tikzcd}[ampersand replacement=\&]
\bullet \arrow[bend left=50]{r}[name=U,below]{} \arrow[bend right=50]{r}[name=D]{} \arrow[Rightarrow,to path=(U) -- (D)]{} \& \bullet \arrow{r} \& \bullet \\
\end{tikzcd} & \begin{tikzcd}[ampersand replacement=\&]
\bullet \arrow{r}\&\bullet \arrow[bend left=50]{r}[name=U,below]{} \arrow[bend right=50]{r}[name=D]{} \arrow[Rightarrow,to path=(U) -- (D)]{} \& \bullet \\
\end{tikzcd} & \begin{tikzcd}[ampersand replacement=\&]
\& \bullet \arrow{rd} \arrow[Rightarrow,shorten >= 10pt,shorten <= 10pt]{dd} \& \\
\bullet \arrow{ru} \arrow{rd} \& \& \bullet \\
\& \bullet \arrow{ru} \& \\
\end{tikzcd} \\ 
\underline{D}_1:= \underline{A}_{1,1}\#_0 \underline{P}_1 & \underline{D}_2:=\underline{P}_1\#_0 \underline{A}_{1,1} & \underline{A}_{2,2} \\
\end{array}
\]

They all have for source and target the $1$-polyplex $\underline{P}_2 :=\bullet \rightarrow \bullet \rightarrow \bullet$, hence one can consider the two $3$-plexes:

 \[ \underline{U}: \underline{D}_1 \rightarrow \underline{A}_{2,2}  \text{ and } \underline{V}:\underline{D}_2 \rightarrow \underline{A}_{2,2}. \]

Indeed, any pair of non-identity parallel positive $2$-polyplexes defines a unique positive $3$-plex between them. For example, $\underline{U}$ is the $3$-polygraph:

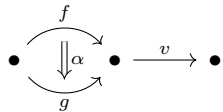
\[\begin{tikzcd}[ampersand replacement=\&]
\bullet \arrow[bend left=50]{r}[name=U]{f} \arrow[bend right=50]{r}[name=D,below]{g} \arrow[Rightarrow,shorten >=5pt,shorten <= 5pt,to path=(U) -- (D)\tikztonodes]{r}{\alpha} \& \bullet \arrow{r}{v} \& \bullet \\
\end{tikzcd} \]

with one additional $2$-cell $\alpha': f \#_0 v \Rightarrow g \#_0 v$ and one additional (universal) $3$-cell $U:\alpha \#_0 v \Rrightarrow \alpha'$.

In particular, one can note that the boundary inclusion $\underline{A}_{2,2} = \pi^+_2 \underline{U} \rightarrow \underline{U}$ is not injective, indeed the two $1$-cells on the right side of $\underline{A}_{2,2}$ have the same image $v$ in $\underline{U}$.

We now form the $3$-polyplex $ \underline{U}\#_0 \underline{A}_{1,1}$ whose underlying polygraph is 

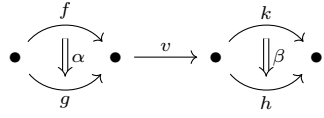
\[\begin{tikzcd}[ampersand replacement=\&]
\bullet \arrow[bend left=50]{r}[name=U]{f} \arrow[bend right=50]{r}[name=D,below]{g} \arrow[Rightarrow,shorten >=5pt,shorten <= 5pt,to path=(U) -- (D)\tikztonodes]{r}{\alpha} \& \bullet \arrow{r}{v} \& \bullet \arrow[bend left=50]{r}[name=U2]{k} \arrow[bend right=50]{r}[name=D2,below]{h} \arrow[Rightarrow,shorten >=5pt,shorten <= 5pt,to path=(U2) -- (D2)\tikztonodes]{r}{\beta} \& \bullet \\
\end{tikzcd} \]

with the same two additional cells $\alpha'$ and $U$ as above. The $2$-target of $\underline{U} \#_0 \underline{A}_{1,1}$ is $\underline{A}_{2,2} \#_0 \underline{A}_{1,1}$, i.e:

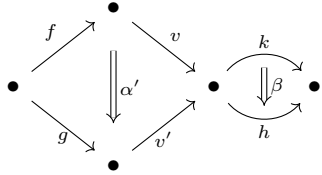
\[\begin{tikzcd}[ampersand replacement=\&]
\& \bullet \arrow{dr}{v} \arrow[Rightarrow,shorten >= 10pt,shorten <= 10pt]{dd}{\alpha'} \& \& \\
\bullet \arrow{ru}{f} \arrow{rd}[below]{g}  \& \& \bullet \arrow[bend left=50]{r}[name=U2]{k} \arrow[bend right=50]{r}[name=D2,below]{h} \arrow[Rightarrow,shorten >=5pt,shorten <= 5pt,to path=(U2) -- (D2)\tikztonodes]{r}{\beta} \& \bullet \\
\& \bullet \arrow{ur}[below]{v'}\& \\
\end{tikzcd} \]

and it is sent to $\underline{U} \#_0 \underline{A}_{1,1}$ by sending each cell in the diagram above to the cell with the same name (with both $v$ and $v'$ are sent to $v$), here again, the map $\pi^+_2(\underline{U} \#_0 \underline{A}_{1,1}) \rightarrow \underline{U} \#_0 \underline{A}_{1,1}$ is not a monomorphism. Our example is obtained by post-composing this $3$-polyplex with (a whiskering of) $V$, but there is two different way we can do that:

We can either apply $V$ on the sub-diagram corresponding to $v,h,k$ and $\beta$, which means forming the composite:

 \[ (\underline{U} \#_0 \underline{A}_{1,1}) \#_2 [(\underline{f} \#_0 \underline{V}) \#_1 (\underline{A}_{2,2} \#_0 \underline{h} )],\]

or apply $V$ to the sub-diagram corresponding to $v',h,k$ and $\beta$, which means forming the composite:

\[( \underline{U} \#_0 \underline{A}_{1,1} ) \#_2 [(\underline{A}_{2,2} \#_0 \underline{k}) \#_1 (\underline{g} \#_0 \underline{V})]. \]

In both these expressions, all the $\underline{f},\underline{h},\underline{g}$ and $\underline{k}$ should all be just ``$\underline{P}_1$'' but we thought that having the name of the precise arrow it corresponds in the diagrams would help the reader to parse the notation. As $v$ and $v'$ are actually collapsed together in the underlying polygraph of $\underline{U}$ this makes the underlying polygraphs of the polyplexes we obtain this way identical and equal to:

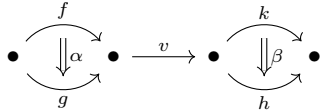
\[\begin{tikzcd}[ampersand replacement=\&]
\bullet \arrow[bend left=50]{r}[name=U]{f} \arrow[bend right=50]{r}[name=D,below]{g} \arrow[Rightarrow,shorten >=5pt,shorten <= 5pt,to path=(U) -- (D)\tikztonodes]{r}{\alpha} \& \bullet \arrow{r}{v} \& \bullet \arrow[bend left=50]{r}[name=U2]{k} \arrow[bend right=50]{r}[name=D2,below]{h} \arrow[Rightarrow,shorten >=5pt,shorten <= 5pt,to path=(U2) -- (D2)\tikztonodes]{r}{\beta} \& \bullet \\
\end{tikzcd} \]

with in addition four $3$-cells:

\[\begin{array}{c c} \alpha': f \#_0 v \Rrightarrow g \#_0 v & \beta': v \#_0 k \Rrightarrow v \#_0 h  \\ 
 U : \alpha \#_0 v \Rrightarrow \alpha' & V : v \#_0 \beta \Rrightarrow \beta' \\
\end{array} \]
 
but they have different universal cells, which corresponds to the two different composite given above, and have different $2$-targets: the first one has for $2$-target $(\alpha' \#_0 k ) \#_1 (g \#_0 \beta')$ while the other has $(f \#_0 \beta') \#_1 (\alpha' \#_0 h)$ and these are different, they even corresponds to different $2$-polyplexes. The point of this story being that $\alpha$ and $\beta$ have disjoint boundary hence the vertical order does not matter when we compose them, but once we replace them with $\alpha'$ and $\beta'$ we still have two ways to vertically compose them but this time their boundaries intersect and hence the vertical composition order matter.

Finally, one also get an example where the a map $\underline{K} \rightarrow \underline{U} \#_2 \underline{K}$ is not injective by simply taking $\underline{K}$ to be the plex whose source and target are both $\underline{A}_{2,2}$, indeed as the map $\underline{A}_{2,2} = \pi_2^+ \underline{ U} \rightarrow \underline{U}$ is not injective while the map $\underline{A}_{2,2} = \pi_2^- \underline{K} \hookrightarrow \underline{K}$ is injective, the gluing of $\underline{K}$ on $\underline{U}$ will force the collapse certain $1$-cell in the image of $\pi_2^- \underline{K}$ (those that are collapsed in $\underline{U}$, i.e. those on the right side of $\underline{A}_{2,2}$) and this makes the map $\underline{K} \rightarrow \underline{U} \#_2 \underline{K}$ non-injective.

}

\block{\label{Ex:SimonForest}After the publication of the first version of this paper, an even more striking example has been found. It is due to Simon Forest (see \cite{forest1018pasting}). It was introduced as a counter-example to some claims in older works on pasting diagrams, but it also gives an example of a positive polygraph which is a positive polyplex in two different way, but with same source and targets arrow in both case. This example is important for two reasons: First, one could believe after the previous example that specifying the underlying polygraph and the source and/or target is sufficient to characterize a polyplex, this would salvage some sort of order independence of the composition (informally, the composition would be order independent as soon as one specifies its source and target as well) but this is not the case. Secondly, it has been suggested several time that polyplexes could be defined as certain higher co-spans in the category of polygraphs (for example, in A.Burroni's texts \cite{burroni2012automates} on the question where such higher co-spans appears under the name ``transiteurs'' or ``logographs''). Indeed given a polyplex $\underline{p}$ the $\pi^{\epsilon}_k \underline{p}$ and maps between them form such a higher co-span. In fact our proof in section \ref{subsection_AutomPD-SPpolyg} actually shows that given a source positive polyplex the corresponding higher co-span of polygraphs does not have any automorphisms. But the following example shows that two different polyplexes can still have isomorphic underlying higher co-spans of polygraphs. 

\bigskip

One starts with the following $2$-polyplex, which will be both the source and target of the two $3$-polyplexes that we will construct:

\[ \underline{G} = (\underline{A} \#_1 \underline{C}) \#_0 ( \underline{B} \#_1 \underline{D} ) :  \parbox{3.8cm}{\begin{tikzcd}[ampersand replacement=\&,column sep=large]
\bullet \arrow[bend left=75]{r}[name=UA,below]{} \arrow[r, " "{name=DA}, " "{name=UC,below}] \arrow[bend right=75]{r}[name=DC]{} \arrow[Rightarrow,to path=(UA) -- (DA)\tikztonodes]{}{A} \arrow[Rightarrow,to path=(UC) -- (DC)\tikztonodes]{}{C} \& \bullet \arrow[bend left=75]{r}[name=UB,below]{} \arrow[r, " "{name=DB}, " "{name=UD,below}] \arrow[bend right=75]{r}[name=DD]{} \arrow[Rightarrow,to path=(UB) -- (DB)\tikztonodes]{}{B} \arrow[Rightarrow,to path=(UD) -- (DD)\tikztonodes]{}{D} \& \bullet \\
\end{tikzcd}} \]
(where $\underline{A},\underline{B},\underline{C},\underline{D}$ are all the plexes $\underline{A}_{1,1}$

One then consider two $3$-plexes:

\[   \underline{U} : \left( \parbox{3.9cm}{\begin{tikzcd}[ampersand replacement=\&,column sep=large]
\bullet \arrow[bend left=75]{r}[name=UA,below]{} \arrow[r, " "{name=DA}, " "{name=UC,below}] \arrow[Rightarrow,to path=(UA) -- (DA)\tikztonodes]{}{A} \& \bullet \arrow[r, " "{name=DB}, " "{name=UD,below}] \arrow[bend right=75]{r}[name=DD]{} \arrow[Rightarrow,to path=(UD) -- (DD)\tikztonodes]{}{D} \& \bullet \\
\end{tikzcd}} \right) \quad  \overset{U}{\Rrightarrow} \quad \left( \parbox{3.9cm}{\begin{tikzcd}[ampersand replacement=\&,column sep=large]
\bullet \arrow[bend left=75]{r}[name=UA,below]{} \arrow[r, " "{name=DA}, " "{name=UC,below}] \arrow[Rightarrow,to path=(UA) -- (DA)\tikztonodes]{}{A'} \& \bullet \arrow[r, " "{name=DB}, " "{name=UD,below}] \arrow[bend right=75]{r}[name=DD]{} \arrow[Rightarrow,to path=(UD) -- (DD)\tikztonodes]{}{D'} \& \bullet \\
\end{tikzcd}} \right) \]

\[ \underline{V} :  \left( \parbox{3.8cm}{\begin{tikzcd}[ampersand replacement=\&,column sep=large]
\bullet \arrow[r, " "{name=DA}, " "{name=UC,below}] \arrow[bend right=75]{r}[name=DC]{} \arrow[Rightarrow,to path=(UC) -- (DC)\tikztonodes]{}{C} \& \bullet \arrow[bend left=75]{r}[name=UB,below]{} \arrow[r, " "{name=DB}, " "{name=UD,below}] \arrow[Rightarrow,to path=(UB) -- (DB)\tikztonodes]{}{B}  \& \bullet \\
\end{tikzcd}} \right) \quad  \overset{V}{\Rrightarrow} \quad \left( \parbox{3.8cm}{\begin{tikzcd}[ampersand replacement=\&,column sep=large]
\bullet \arrow[r, " "{name=DA}, " "{name=UC,below}] \arrow[bend right=75]{r}[name=DC]{} \arrow[Rightarrow,to path=(UC) -- (DC)\tikztonodes]{}{C'} \& \bullet \arrow[bend left=75]{r}[name=UB,below]{} \arrow[r, " "{name=DB}, " "{name=UD,below}] \arrow[Rightarrow,to path=(UB) -- (DB)\tikztonodes]{}{B'}  \& \bullet \\
\end{tikzcd}} \right) \]

One can then form two different polyplexes by composing them in different ways:

\begin{itemize}

\item Either one whiskers $\underline{U}$ by $\underline{C}$ and $\underline{B}$ so that its source becomes $\underline{G}$, to then postcomposed it by $\underline{V}$ whiskered by $\underline{A'}$ and $\underline{D'}$,

\item Or one whiskers $\underline{V}$ by $\underline{A}$ and $\underline{D}$ so that its source becomes $\underline{G}$, to then postcomposed it by $\underline{U}$ whiskered by $\underline{C'}$ and $\underline{D'}$.

\end{itemize}

In both case the underlying polygraph is the same: one has six $1$-cells and height $2$-cells as follow:

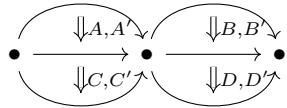
\[ \parbox{3.8cm}{\begin{tikzcd}[ampersand replacement=\&,column sep=large]
\bullet \arrow[bend left=75]{r}[name=UA,below]{} \arrow[r, " "{name=DA}, " "{name=UC,below}] \arrow[bend right=75]{r}[name=DC]{} \arrow[Rightarrow,to path=(UA) -- (DA)\tikztonodes]{}{A,A'} \arrow[Rightarrow,to path=(UC) -- (DC)\tikztonodes]{}{C,C'} \& \bullet \arrow[bend left=75]{r}[name=UB,below]{} \arrow[r, " "{name=DB}, " "{name=UD,below}] \arrow[bend right=75]{r}[name=DD]{} \arrow[Rightarrow,to path=(UB) -- (DB)\tikztonodes]{}{B,B'} \arrow[Rightarrow,to path=(UD) -- (DD)\tikztonodes]{}{D,D'} \& \bullet \\
\end{tikzcd}}\]

And two $3$-cells given by $U$ and $V$ with source and target as specified above. The two ways of composing $U$ and $V$, corresponding to the two polyplexes above, are distinct arrows of this polygraph (the sources and targets of $U$ and $V$ are not disjoint so there is no way to exchange their composition using some sort of Eckmann-Hilton argument, an actual proof that they are distinct will be available in \cite{forest1018pasting}) and hence one indeed have two different polyplexes, but in both case their source corresponds to the image of $\underline{G}$ sent on $A,B,C,D$ and their target is the image of $\underline{G}$ sent on $A',B',C',D'$. In particular they both have the same underlying higher co-span.
}

\block{\label{Discuss_Rel_Batanin}We would like to finish this section with a rather technical observation for the expert reader that play no concrete role in the paper: We want to show that M.Batanin criterion in \cite{batanin2002computads} for proving that the category of computads associated to a globular operad is a presheaf category is insufficient for the main case of interest to us, and try to explain the relation between our theorem and this criterion. We refer the reader to \cite{batanin2002computads} for the notions that we will mention below.

We start with a quick remark: in his paper M.Batanin quote a result of A.Carboni and P.T.Johnstone from \cite{carboni1995connected} saying that a finitary monad on the category of sets is familially representable if and only if it corresponds to a strongly regular theory. This result is unfortunately the other false result of \cite{carboni1995connected} (see \cite{carboni2004corrigenda} for a counterexample): instead, strongly regular theories correspond to non-symmetric operads while the condition of familial representability corresponds to the weaker notion of $\Sigma$-cofibrant symmetric operad. Due to the strong similarities between M.Batanin results and the results of the present paper, we strongly believe\footnote{We have unfortunately not being able to understand M.Batanin proofs well enough to be sure of that, but he does not seem to be using strong regularity anywhere.} that all his results in \cite{batanin2002computads} are correct if we replace everywhere in his paper ''strongly regular'' by this weaker condition of being a familially representable monad or equivalently being a $\Sigma$-cofibrant operad.

\bigskip

This being said, the class of positive polygraphs can be seen as the category of polygraphs for the globular operad for non-unital $\infty$-categories (see \ref{Def_nonUnitalInfinityCat} for our notion of non-unital $\infty$-category) which is just the sub-monad of the free $\infty$-category monad $D$ on the category of globular sets defined by 

\[ M(X)=\{a \in D(X)| a \text{ is not an identity arrow in D(X).}\} \]

One easily check from an explicit description of $D$ (See for example \cite[III.8.1]{leinster2004higher}) that $M(X)$ defined this way is indeed a sub-globular set of $D(X)$ and that $M$ is a sub-monad of $D$, with the inclusion $M \rightarrow D$ being cartesian. Hence $M$ is a globular operad in the sense defined in \cite{batanin2002computads}, and the $M$-computads in his sense are exactly our positive polygraphs. But the second slice $\mathcal{P}_2(M)$ of $M$ is not a $\Sigma$-cofibrant operad. Indeed a $\mathcal{P}_2(M)$-algebra is the same as a non-unital $2$-category (in the sense of an $M$-algebra, or in the sense of \ref{Def_nonUnitalInfinityCat}) with only one cell in dimension $0$ and dimension $1$, but this is enough to form a kind of ``Eckmann-Hilton collapse'': if $x$ and $y$ are two $2$-arrows of such an $\infty$-category and let $e$ denotes the unique $1$-arrow of this category then, as $e$ is a $1$-arrow one has $e \#_1 x= x \#_1 e =x$ and as there is only one $1$-arrow one has $e \#_0 e =e$, from there:

\[\begin{array}{c c c} (e \#_0 x \#_0 e) \#_1 (e \#_0 y \#_0 e) & = & [(e \#_0 x) \#_0 e] \#_1 [e \#_0 (y \#_0 e)] \\
& = & [(e \#_0 x) \#_1 e] \#_0 [e \#_1 ( y \#_0 e)]  \\
& = & e \#_0 x \#_0 y \#_0 e \\
\end{array} \]

and starting from the other possible bracketing one gets:

\[\begin{array}{c c c} (e \#_0 x \#_0 e) \#_1 (e \#_0 y \#_0 e) & = & [e \#_0 (x \#_0 e)] \#_1 [(e \#_0 y) \#_0 e] \\
& = & [e \#_1 (e \#_0 y)] \#_0 [(x \#_0 e) \#_1 e ]  \\
& = & e \#_0 y \#_0 x \#_0 e \\
\end{array} \]

hence one has a commutative operation on two variables:

\[ e \#_0 y \#_0 x \#_0 e = e \#_0 x \#_0 y \#_0 e \]

which shows that the corresponding operad is not $\Sigma$-cofibrant: it admits a generic\footnote{By generic operations, we mean an operation where each variables that appears in it appears exactly once. For an ordinary finitary algebraic theory, this only makes sense if the theory corresponds to a symmetric operads (which is the case of such ``slice'' theory by results of M.Batanin) and in this case the generic operations on $n$-variables are exactly the operation in the set $\Ocal(n)$ of the operads.} operation invariant under a non trivial permutation of its variables.

Informally, what happens here is that in M.Batanin criterion we look at how the $n$-dimensional operations which are generic on variables of dimension $n$ but where all the lower dimensional variables are set to be equals behave under permutations of the variables of dimension $n$, while in our framework we are looking at the behavior under permutation of the variables of operations which are ``globally generic'' i.e. where in some sense all the variables of all dimension are used exactly once (this is a very vague formulation, which is made rigorous by the notion of polyplex). In both case the criterion is that in the absence of such permutations one has a presheaf category. Moreover, we prove in our context (lemma \ref{Prop_AutActsOnTopCellsOnly}) that it is enough to show that such permutations acts trivially on the $n$-dimensional arrows, which allows to recover M.Batanin criterion for sub-operads of $D$. But on the other hand, our proof in sub-section \ref{subsection_AutomPD-SPpolyg} relies heavily one carefully analyzing how such permutations of a polyplex acts on $(n-1)$-variables (i.e. $(n-1)$-cells of the polyplexes), which are all collapsed together in the framework of M.Batanin criterion.

\bigskip

But obviously, while our criterion is more powerful in the sense that it is a necessary and sufficient condition at least for detecting good classes of polygraphs, M.Batanin criterion applies to the more general situation of any globular operads. As mentioned earlier, we hope to extend our main theorem \ref{Th_mainGoodClass} to any globular operads, or even to a more general notion of operad using other non-globular sort of combinatorics (typically to any parametric right adjoint cartesian monad on a category of presheaves over a directed category).
}

\appendix

\renewcommand{\thesubsubsection}{\Alph{section}.\arabic{subsubsection}}

\section{On C.Simpson's conjecture and the Kapranov-Voevodsky strategy}
\label{App-Simpson-voevodsky-Kapranov}

\blockn{In 1991, M.Kapranov and V.Voevodsky published a proof (\cite{kapranov1991infty}) of a form of the homotopy hypothesis, claiming that the homotopy category of spaces is equivalent to the homotopy category of a certain kind of $\infty$-groupoids, which are strict $\infty$-categories where every arrow is weakly invertible in every degree (we refer to \cite{kapranov1991infty} for the precise definition of weak invertibility). In 1998, C.Simpson published a proof (\cite{simpson1998homotopy}) that this cannot be true. He was not able to point out a precise mistake in \cite{kapranov1991infty} but he conjectured that this had to do with how units are handled and he formulated the conjecture now called the Simpson conjecture or Simpson's semi-strictification conjecture, that the homotopy category of spaces is equivalent to the homotopy category of a notion of $\infty$-groupoids where the associativity and the exchange rule are strict, but units and inverses are weak. He did not gave a precise definition of what this means, and he added that he ``\emph{thinks that the argument of \cite{kapranov1991infty} [...] actually serves to prove the above statement}''.

At the time the present paper is written it seems to be largely accepted that it is C.Simpson's paper which is correct and M.Kapranov and V.Voevodsky's paper which is flawed, but as far as we know it is still unclear where is the mistake in this paper. Also C.Simpson's conjecture is still open\footnote{As far as we know, and with the exception of \cite{henry2018regular} which will be mentioned latter, the only concrete progress since its formulation is a proof of a form of the conjecture for $3$-groupoid with only one object in \cite{joyal2006weak}.}, and while it is still plausible that the general strategy of \cite{kapranov1991infty} could gives a proof, it seems that it needs an in depth reworking and some new ideas in order to achieve that. The results about positive polygraphs that we proved in the present paper somehow steam from my personal analysis of ``why'' the paper \cite{kapranov1991infty} fails and what should be done to fix it in order to prove Simpson's conjecture.

The goal of this section is to sum up this analysis, both in order to motivate the present paper and to explain our plan to attack this conjecture. A subsequent paper (\cite{henry2018regular}) will make the ideas presented here more precise, and push some of these further to prove a form of Simpson's conjecture in all dimension. There are still some other form of the conjecture that are open.
}

\block{We start by explaining the general idea of \cite{kapranov1991infty}, which we refer to as the ``Kapranov-Voevodsky strategy''. One starts with a topological space $X$ and we want to define a ``fundamental $\infty$-groupoid of $X$'', whose objects are points of $X$, $1$-arrows are paths between points, $2$-arrows are endpoints preserving homotopies, and more generally $n$-morphisms are boundary preserving homotopies between $(n-1)$-morphisms. In this groupoid compositions are given by compositions of homotopies, but we want to define it in such a way that associativity and exchange rule holds strictly, and not up to higher homotopies as it would be the case with a naive definition.

\bigskip

If we just look at $1$-arrows, composition can be made strict by using ``Moore paths'' i.e. by allowing the length of paths to vary: if one has two composable $1$-arrows given by $[0,n] \rightarrow X$ and  $[0,m] \rightarrow X$ and we define their composite as a map $[0,n+m] \rightarrow X$ then one gets a strictly associative composition. The starting idea in \cite{kapranov1991infty} is to push this idea to all dimensions: $2$-arrows now need to be homotopies between paths of possibly different length and if we want to compose them in a strictly associative way (and with strict exchange rules) the only way is to define the composition formally by just gluing together the spaces indexing those homotopies. To generalize this in all dimensions, we need to come up with a notion of ``generalized Moore homotopy'' indexed by certain ``diagrams'' or ``cell complexes''. The introduction of \cite{kapranov1991infty} give a nice explanation of this idea.

\bigskip

M.Kapranov and V.Voevodsky propose to use M.Johnson's notion of pasting diagram from \cite{johnson1989combinatorics} which they improved in a companion paper \cite{kapranov1991combinatorial}. They define a category of Johnson's diagrams and a geometric realization functor for such diagrams. One can then attempt to make the above idea formal by saying that we want the $n$-arrows of the fundamental $\infty$-groupoid of $X$ to be pairs of a pasting diagram $K$ together with a continuous map $|K| \rightarrow X$ where $|K|$ denotes the geometric realization of the diagram $K$, and composition should be defined by pasting of diagrams. Unfortunately such a construction have no chance of being either a left of right adjoint functor, and proving that a functor induces an equivalence between homotopy categories can be very hard if this functor is not part of an adjunction.

To avoid this problem, as well as possibly some technical difficulties in making the above definition formal, they move to a slightly different construction: if $D$ denotes their category of Johnson diagrams, then they have two functors $D \rightarrow \infty-cat$ and $D \rightarrow Spaces$ respectively defined as the free $\infty$-category generated by a diagram and the geometric realization functor. Using usual Kan extension techniques one gets two left adjoint functors $Prsh(D) \rightarrow \infty-cat$ and $Prsh(D) \rightarrow Spaces$. They then use the fact that $Prsh(D)$ is somehow similar to the category of simplicial sets or cubical sets to set up a homotopy theory on this presheaf category, and they claim to prove that these two adjunctions induce equivalences between the three homotopy categories of interest, which concludes their proof.
}

\block{If we forget the more indirect version of the construction they actually uses, and come back to the initial idea of using some ``generalized Moore homotopies'' indexed by some class of pasting diagrams, then one can see that in order to define a well behaved fundamental $\infty$-groupoid of $X$ whose $n$-arrows are maps $|K| \rightarrow X$ for $K$ a pasting diagram, it seem to us that one should at least expect that the category of diagrams we are using satisfies the following properties:

\begin{itemize}

\item[(A)] One should be able to ``compose'' the pasting diagrams that we use, by gluing them together, and geometric realization should be compatible to these gluing. We need this to define the $\infty$-categorical composition operations on the set of maps $|K| \rightarrow X$.

\item[(B)] If $K$ and $K'$ are two ``parallel'' $n$-dimensional pasting diagrams (their boundaries are the same diagrams) then one should be able to construct a new pasting diagram by gluing $K$ and $K'$ together along their boundary and adding one new $(n+1)$-cell between them. Having this, allows us to comfortably see inside this ``fundamental $\infty$-groupoid of $X$'' that if two parallel cells corresponding to maps $|K| \rightarrow X$ and $|K'| \rightarrow X$ are homotopic in $X$ then this homotopy is detected inside the fundamental $\infty$-groupoid.

\end{itemize}

\bigskip

It appears that Johnson diagrams fail to have either property $(A)$ or $(B)$.

Property $B$ fails because of the following ``stupid'' example: If one considers the diagram $\bullet \rightarrow \bullet$ representing a path and the diagram $\bullet$ representing the constant path, then a $2$-arrow between them would be diagram with a loop and a contraction of that loop. But Johnson's diagrams are not allowed to contains loop (they satisfies a certain condition called ``loop free'' which as the name suggest in particular implies that the underlying $1$-graph cannot have loops). This first obstruction is clearly related to the presence of units and disappear when we work in a ``non-unital framework'' as suggested by C.Simpson.

While this first observation is very encouraging for C.Simpson conjecture, there is unfortunately, a second type of counterexample to property $(B)$ that still exists even when we restrict ourselves to non-identity arrows. Consider the following two Johnson $2$-pasting diagrams:

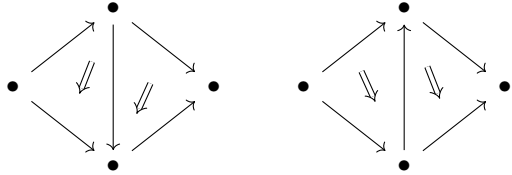
\[\begin{tikzcd}[ampersand replacement=\&]
\&  |[alias=T]| \bullet \arrow{dd} \arrow{dr}[name=S]{} \& \\
\bullet \arrow{ur} \arrow{dr}[name=D]{} \& \& \bullet \\
\& |[alias=TT]| \bullet \arrow{ur} \& \arrow[Rightarrow,from= T, to =D,shorten >= 8pt,shorten <= 15pt]{}  \arrow[Rightarrow,from= S, to =TT,shorten >= 15pt,shorten <= 15pt]{} \\
\end{tikzcd} \qquad \begin{tikzcd}[ampersand replacement=\&]
\&|[alias=S]| \bullet \arrow{dr} \& \\
\bullet \arrow{ur}[name=SA]{} \arrow{dr} \& \& \bullet \\
\& |[alias=T]| \bullet \arrow{ur}[name=TA]{} \arrow{uu} \& \arrow[Rightarrow,from= SA, to =T,shorten >= 20pt,shorten <= 10pt]{} \arrow[Rightarrow,from=S, to =TA,shorten >= 6pt,shorten <= 17pt]{} \\
\end{tikzcd} \]

They are both legitimate Johnson diagrams and they are parallel, but if we glue them along their common boundary, the two vertical arrows will again form a loop, so this takes us outside of the class of Johnson diagrams. More problematically, for the exact same reason there can be no Johnson $3$-diagrams which has these two diagrams as source and target. So if one form a fundamental $\infty$-groupoid of $X$ whose arrows are only parametrized by Johnson diagrams there can never be a $3$-arrow between two $2$-arrows parametrized by the two pasting diagram above. Hence this $\infty$-groupoid will always ``miss'' certain higher homotopies between $n$-arrows.

Moreover, this also provides a counterexample for condition $(A)$, one can consider the following two $3$-dimensional Johnson's diagram:

\begin{itemize}

\item The diagram representing a single $3$-arrow from the first of these two diagrams to the diagram $\underline{A}_{2,2}$ with the same boundary by only one internal $2$-cell from the source to the target.

\item The diagram representing a single $3$-arrow from $\underline{A}_{2,2}$ to the second of these two diagrams.  
\end{itemize}

These two diagrams are indeed Johnson pasting diagrams, but if we compose them (along $\underline{A}_{2,2}$) one again gets a loop with the two vertical arrows, so the composite is not a Johnson diagram ! Hence it is not even possible to define an $\infty$-category this way.

Although not explicitly claimed in their paper, M.Kapranov and V.Voevodsky seem to use that Johnson's pasting diagrams can be composed, at least in the proof of their lemma $3.4$ of \cite{kapranov1991infty}. This might be one of the problems with their proof.

 }

\block{We believe that at the end of the day, and even if that might not be the exact technical reason which makes their proof incorrect, the reason why the original form of the Kapranov-Voevodsky strategy cannot succeed is exactly because of the failure of these two properties for Johnson diagrams, and that fixing the proof in order to proves C.Simpson's conjecture requires to construct a new class of diagrams which satisfies those two properties.

\bigskip

What the present paper achieve is exactly to construct such a category of diagrams in the ``non-unital'' case. Indeed, positive polyplexes are exactly the class of diagram that is ``generated'' by these construction $(A)$ and $(B)$ above, in the non-unital case. We do not know how to do that in the unital case, and it is likely to be impossible (it is indeed provably impossible as soon as we impose more precise conditions on how those diagrams should behave). So as a first step one should replace the category of Johnson diagrams with the category of ``positive polyplexes'' in the sense of the present paper, but we will see below that there is something even more natural to do: instead of looking at presheaves on the category of polyplexes, we will look at presheaves on the category of plexes, i.e. just the category of positive polygraphs.}

\block{One of the problem with C.Simpson's conjecture is that the notion of $\infty$-category with weak units, or of non-unital $\infty$-categories does not have a unique definition. For example, one could understand ``non-unital categories'' as being a globular sets with all the operations $f \#_k g$ defined for $f$ and $g$ of dimension $n$ and $k<n$. Or one can also require to have compositions like $f \#_k g$ defined even when $f$ and $g$ have different dimensions, allowing to define whiskering as $a\#_0 f$:

\[ \begin{tikzcd}[ampersand replacement=\&]
\bullet \arrow[bend left=50]{r}[name=U,below]{} \arrow[bend right=50]{r}[name=D]{} \arrow[Rightarrow,to path=(U) -- (D)\tikztonodes]{r}{\alpha} \& \bullet \arrow{r}{f} \& \bullet \\
\end{tikzcd} \]

instead of using an identity arrow $1_f$ to define it as a horizontal composite $a \#_0 1_f$ of two $2$-arrows. One can also use even different composition shapes, for example one can consult \cite{kock2006weak} for a simplicially based definition of weakly unital $\infty$-categories, and all those notions are non-equivalent (in the $1$-categorical sense). So we need to choose such a notion.

Fortunately, the notion of positive polygraphs, for which we proved the existence of a well behaved notion of pasting diagrams (the positive polyplexes), are precisely a notion of polygraph corresponding to a certain well defined notion of non-unital $\infty$-category: Globular sets where all the compositions $f \#_k h $ are defined even when $f$ and $g$ have different dimensions, and all associativities and exchange rules holds. Almost surprisingly, when writing down the definition of such ``non-unital $\infty$-categories'' one sees that they are in fact ordinary $\infty$-categories satisfying two additional axioms:

}

\block{\label{Def_nonUnitalInfinityCat}\Def{A non-unital $\infty$-category is an $\infty$-category $X$ in the sense of definition \ref{Def_StreetInfinityCat} which satisfies the following two additional axioms:

\begin{itemize}

\item If $f$ is an arrow of dimension greater than $n$, then the $\pi^{\epsilon}_n(f)$ are of dimension $n$ exactly.

\item for any $f,g,k$ such that $f \#_k g$ is defined, its dimension is the maximum of the dimension of $f$ and the dimension of $g$.

\end{itemize}

A morphism of non-unital $\infty$-categories is a dimension preserving morphism of $\infty$-categories. The category of non-unital $\infty$-categories is denoted $\infty-Cat^{nu}$.
}

Of course the interpretation of the set $X$ has changed: it is no longer the increasing union of the set of $n$-arrows, but the \emph{disjoint union} of the sets of $n$-arrows. The fact that a non-unital $\infty$-category in this sense can be identified with a certain ordinary $\infty$-category corresponds just to the fact that there is a faithful ``unitarization functor'' from $\infty-cat^{nu}$ to $\infty-cat$ which just ``freely add units''.

We claim that those $\infty$-categories are the algebras for a monad on globular sets (in fact the monad $M$ mentioned in \ref{Discuss_Rel_Batanin}), and that the polygraphs for this monad (following \cite{batanin1998computads}) are exactly the positive polygraphs. In particular, one has an adjunction $( \_ )^*: \Pb^+ \leftrightarrows \infty-cat^{nu} : N$, and $(\_)^*$ is just the usual free $\infty$-category on a polygraph, which happens to takes values in $\infty-cat^{nu}$ when applied in the category of positive polygraphs.

\bigskip

Nonetheless, the fact that those non-unital categories are a (non-full) subcategory of the category of strict $\infty$-categories means that, if the form of C.Simpson's that we are going to conjecture holds, then the main result of M.Kapranov and V.Voevodsky was a lot closer to be true than what we thought: every homotopy type would be representable by a strict unital $\infty$-category: the unitarization of its non-unital $\infty$-groupoids, and this $\infty$-category indeed computes the correct homotopy groups, we just need to use a definition of homotopy groups that do not use the canonical identity arrows of the category, but a ``weak'' units instead.
}

\block{Finally, again following the footstep of Kapranov and Voevodsky, one defines what we will call the ``Naive geometric realization functor'': $\Pb^+ \rightarrow Spaces$ which send any polygraph $P$ to the geometric realization of the category $Plex^+/P$ of cells of $P$. This is a left adjoint functor and one has a diagram of left adjoint functor:

\[\begin{tikzcd}[ampersand replacement=\&]
Spaces \& \& \infty-cat^{nu} \\
\& \Pb^+ \arrow{lu}{ | \_ |} \arrow{ru}[below]{(\_)^*} \& \\
\end{tikzcd}\]

With $\Pb^+$ being a presheaf category as proved in the present paper.

This could be our new basis to make the Kapranov-Voevodsky strategy into a proof of Simpson's conjecture, and it has a very nice new feature that is not present in the original Kapranov-Voevodsky strategy: if one start with a topological space, applies to it the right adjoint functor from Spaces to $\Pb^+$ and then the free $\infty$-category functor, we obtain exactly the $\infty$-category corresponding to the intuitive idea of generalized Moore spaces we started from:

Indeed if we denote by $N$ the right adjoint to the geometric realization, a cell of $N(X)^*$ is a map from a polyplex $\underline{p}$ to $N(X)$ which is exactly the same as a continuous map from $|\underline{p}|$ to $X$. Hence cells of the corresponding $\infty$-category are exactly maps from the geometric realization of a pasting diagram (a positive polyplex) to the space $X$ as expected in the beginning.}

\block{But unfortunately there is a new problem that comes with the increased complexity in the shapes of the diagrams that we use: Basically, the ``naive'' geometric realization is too naive, and cannot be used in this pictures. To clarify the following discussion we will admit the following, which are proved in the subsequent paper \cite{henry2018regular}:

\begin{itemize}

\item One can construct a ``weak model structure'' (in the sense of \cite{henry2018weakmodel}) on $\infty-cat^{nu}$ where the fibrant objects are non-unital $\infty$-categories in which every arrow has a weak identity endomorphism (which are defined as weakly idempotent endomorphisms) and which satisfies all of Kapranov and Voevodsky divisibility condition of \cite{kapranov1991infty}. The weak equivalences between fibrant objects being the map inducing a bijection on all the $\pi_n$.

\item One can construct a similar ``weak model structure'' on $\Pb^+$ where the fibrant objects are polygraphs $P$ such that $P^*$ is fibrant in the previous sense and such that for each arrow $f \in P^*$ there exists a cell $f' \in P$ parallel to $f$ and a cell $a$ between $f$ and $f'$. Weak equivalences in $\Pb^+$ are the arrows that induces equivalences in the sense of the above weak model structure on $\infty-cat^{nu}$

\item The pair of adjoint functor between $\Pb^+$ and $\infty-cat^{nu}$ defined above is a Quillen equivalence (for the notion of Quillen equivalence adapted to weak model structures introduced in \cite{henry2018weakmodel}). 

\end{itemize}

We refer to \cite{henry2018regular} for proof and more precise statement of these claims.
}

\block{We can now explain the problem that remains to be solved in order to prove Simpson's conjecture: the naive geometric realization is not even a Quillen functor. The reason for this is that the naive geometric realization send any plex to a contractible topological space. But quite surprisingly, and contrary to what happen with Johnson diagrams in \cite{kapranov1991infty}, general plexes are not at all nice contractible balls ! The following is an example of \emph{a non-contractible $3$-plex} : Consider first the following $2$-polyplex:

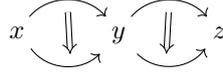
\[ \begin{tikzcd}[ampersand replacement=\&]
x \arrow[bend left=50]{r}[name=U,below]{}
\arrow[bend right=50]{r}[name=D]{} \&
y \arrow[Rightarrow,to path=(U) -- (D)]{}  \arrow[bend left=50]{r}[name=V,below]{}
\arrow[bend right=50]{r}[name=W]{} \& z \arrow[Rightarrow,to path=(V) -- (W)]{}
\end{tikzcd} \]
 
One can form its unique $3$-plex ``endomorphism'', which has the following underlying polygraph:

\begin{itemize}

\item It has three $0$-cells $x,y,z$.

\item It has four $1$-cells: $f,g:x \rightrightarrows y$ and $h,k:y \rightrightarrows z$.

\item It has four $2$-cells: $\alpha,\beta:f \rightrightarrows g$, $\gamma,\delta:h \rightrightarrows k$.

\item It has one $3$-cells: $\Omega: \alpha \#_0 \gamma \rightarrow \beta \#_0 \delta $

\end{itemize}

we claim that this plex is not homotopy equivalent to a point in the model structure we mentioned above. 

The reason for that is that removing a given $n$-cell $h$ and an $(n-1)$-cell $a$ appearing exactly once in the source or the target of $h$ (and not in the other) does not change the homotopy type of a polygraph in this weak model structure (the reader can note that those corresponds exactly to how the ``generating trivial cofibration'' of \cite{kapranov1991infty} are defined in the framework of Johnson's diagram). Admitting that, one can gradually remove the following pairs of cells to the plex above without changing its homotopy type: $(\Omega,\gamma)$ , $(\delta,k)$ $(h,z)$, after that it only remains the following cells: two $0$-cells: $x,y$; two $1$-cells $f,g:x\rightrightarrows y$ and two $2$-cells $\alpha,\beta:h \rightrightarrows k$, i.e. it is exactly the globular polygraph corresponding to a free pair of parallel $2$-cells, whose geometric realization is the $2$-sphere, and is not contractible. 

Obviously this makes no real sense unless we introduce the details of this model structure as well, but we can also see that this polygraph should have the homotopy type of a $2$-sphere in a less formal but more intuitive way: if we think of it in topological terms, its boundary corresponds to two $2$-spheres glued together on a point, or equivalently a single ``twice bigger'' $2$-sphere whose equator has been contracted to a single point. Adding the unique $3$-dimension cell should corresponds to gluing a $3$-ball on this $2$-sphere, i.e. to fill the interior of that single sphere. One hence obtains a $3$-ball whose equator is contracted to a single point. But only the boundary of the equatorial disk is contracted, not the whole $2$-dimensional disk. The resulting space can be deformed into its equatorial disk, which is a two dimensional disk whose boundary has been contracted to point, i.e. a $2$-sphere as claimed above !

The fact that we are able to understand this homotopy type both topologically and in terms of our weak model structure suggest that this problem is not an obstruction to Simpson's conjecture, but only a sign that we need a more subtle geometric realization functor, we hence propose the following conjectures, which basically form a more precise version of C.Simpson conjecture: }

\block{

\Conjecture{There exists a geometric realization functor $| \_ |$ from $\Pb^+$ to Spaces (topological spaces, or simplicial sets) such that:

\begin{itemize}

\item The geometric realization of the polygraph $P_0$ (with just a single cell) is a point.

\item $| \_ |$ is a left adjoint functor.

\item $| \_ |$ send monomorphisms in $\Pb^+$ to cofibrations of spaces.

\item If $P$ is an $n$-plex for $n>0$ and $A$ is the sub-polygraph of $P$ obtained by removing the unique $n$-dimensional cell of $P$ and a single $(n-1)$-dimensional cell of $P$, then the cofibration $| A | \hookrightarrow |P|$ is a weak equivalence.

\end{itemize}

}

Such a functor would give a left Quillen functor from $\Pb^+$ to Spaces sending the point to the point. Hence if our version of the Simpson conjecture holds, and $\Pb^+$ is indeed Quillen equivalent to the model category of spaces, then the usual ``universal property'' of the model $\infty$-category of spaces should implies that:

\Conjecture{Any functor satisfying the condition of the first conjecture is a left Quillen equivalence.}

And these two conjectures, (together with the claim that we made earlier on the existence of weak model structures), implies C.Simpson's conjecture. We believe the hard part is the first conjecture, i.e. constructing a good geometric realization functor}

\block{Another possible approach to completely circumvent this difficulty is to restrict the class of polygraphs we are using to avoid these non-contractible plexes. This is what we do in \cite{henry2018regular}. This corresponds to restricting the type of compositions allowed in our non-unital categories to diagrams which are ``topologically balls''. We call such compositions ``regular'', and the corresponding notion of $\infty$-category ``regular $\infty$-category'' (which are hence not quite $\infty$-categories as they have less composition operations defined), the precise definition of these notions being in \cite{henry2018regular}. Things like whiskerings and horizontal compositions:

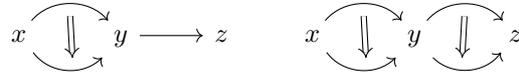
\[ \begin{tikzcd}[ampersand replacement=\&]
x \arrow[bend left=50]{r}[name=U,below]{}
\arrow[bend right=50]{r}[name=D]{} \&
y \arrow[Rightarrow,to path=(U) -- (D)]{} \arrow{r} \& z
\end{tikzcd} \qquad \begin{tikzcd}[ampersand replacement=\&]
x \arrow[bend left=50]{r}[name=U,below]{}
\arrow[bend right=50]{r}[name=D]{} \&
y \arrow[Rightarrow,to path=(U) -- (D)]{}  \arrow[bend left=50]{r}[name=V,below]{}
\arrow[bend right=50]{r}[name=W]{} \& z \arrow[Rightarrow,to path=(V) -- (W)]{}
\end{tikzcd} \]

are not regular, and would \emph{not} be defined in a regular $\infty$-category, but compositions of more complex shapes, possibly containing them, like:

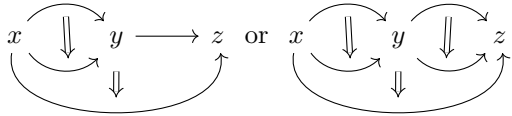
\[ \begin{tikzcd}[ampersand replacement=\&]
x \arrow[bend right=100]{rr}[name=k]{} \arrow[bend left=50]{r}[name=U,below]{}
\arrow[bend right=50]{r}[name=D]{} \&
y \arrow[Rightarrow,to path=(U) -- (D)]{} \arrow[Rightarrow,from=1-2,to=k,shorten >= 1pt,shorten <= 5pt]{}  \arrow{r} \& z 
\end{tikzcd} \text{ or }  \begin{tikzcd}[ampersand replacement=\&]
x \arrow[bend right=100]{rr}[name=k]{} \arrow[bend left=50]{r}[name=U,below]{}
\arrow[bend right=50]{r}[name=D]{} \&
y \arrow[Rightarrow,to path=(U) -- (D)]{} \arrow[Rightarrow,from=1-2,to=k,shorten >= 1pt,shorten <= 5pt]{}  \arrow[bend left=50]{r}[name=V,below]{}
\arrow[bend right=50]{r}[name=W]{} \& z \arrow[Rightarrow,to path=(V) -- (W)]{}
\end{tikzcd} \]

are regular, and would be defined, and ``associative'' in the sense that any two way of composing such diagrams would be equals. As soon as one has weak units (still defined as weakly idempotent endomorphisms) One can recover weak form of whiskering and horizontal composition by inserting weak units in such diagrams.

The form of Simpson's conjecture for this type of $\infty$-category is in some sense weaker and appears easier, and we will prove it in \cite{henry2018regular}.

}

\block{Finally, we clarify the main three differences between our approach and the original approach of \cite{kapranov1991infty}:

\begin{itemize}

\item One has changed the category of diagrams.

\item We use the category of presheaf on $Plex^+$ (i.e. the composition diagrams with just a single $(n+1)$-cell between two $n$-pasting diagrams) instead of the category of all pasting diagrams (polyplexes).

\item We are not introducing degeneracies in our category of diagrams.

\end{itemize}

We already discussed in length why the change in the category of diagrams was necessary, and we explained that there seems to be an incorrect assumption used in lemma $3.4$ of \cite{kapranov1991infty} regarding the fact that Johnson's diagrams can be composed. Moreover as the goal is to remove units, not introducing degeneracies seems very natural as those essentially corresponds to unit. One should also mention that the absence of degeneracies is the reason why we need to move to weak model structures instead of Quillen model structures: for example it is well known that there is no model structure on semi-simplicial sets (simplicial sets without degeneracies) where the (trivial) fibrations are the usual Kan (trivial) fibrations and the weak equivalences are the homotopy equivalences, but we will show in \cite{henry2018weakmodel} that such a weak model structure can be constructed on semi-simplicial sets.

But one might wonder what is the meaning of replacing presheaves on the category of all pasting diagrams to presheaves on this smaller category. We believe that, if correctly taken into account this is a completely unessential change. But it seems that it has not been correctly taken into account in \cite{kapranov1991infty} and that using the category of all pasting diagrams is actually responsible for at least one direct mistake in their paper: it seems that this makes the lemma $3.4$ already mentioned above trivially false. This lemma claims in particular that if $X$ is a presheaf on their category of Johnson diagrams and $X^*$ denotes the $\infty$-category generated by $X$ (the image of $X$ by the left Kan extension of the natural ``free category functor'' from Johnson diagrams to $\infty$-cat) then any arrow of $X^*$ can be represented by a cell of $X$, i.e. an element of $X(a)$ for $a$ a pasting diagram.

We already mentioned that the the proof of lemma $3.4$ seems to use that Johnson diagram can be composed, which is not the case, but another problem is that, even in situation where all the Johnson pasting diagrams appearing can be composed there is still no way to compose cells of a presheaf on the category of Johnson diagrams, even in the $1$-dimensional situation:

Take the presheaf $C$ obtained as a gluing of two copies of the representable object $\bullet \rightarrow \bullet$ glued along the representable object $\bullet$. The free $\infty$-category obtained is just the $1$-category with two arrow $\bullet \rightarrow \bullet \rightarrow \bullet$. But colimits in presheaf categories are computed objectwise, so if $A$ is any Johnson pasting diagram any map from $A$ to $C$ has to factor in one of the two maps $\bullet \rightarrow \bullet $, so no such maps can ever represent the composite of the two arrows in the free $\infty$-category generated by $C$ !

As far as we know there are at least two solutions to this problem: The first is the one we proposed above to restrict to presheaves on the category of pasting diagram with only one top dimensional cells, and having the other pasting diagrams represented by some gluing of those representable pasting diagrams. In this presheaf category the object $C$ is the same as the pasting diagram $\bullet \rightarrow \bullet \rightarrow \bullet $ and hence the problem disappear. Another alternative would be to restrict to fibrant objects (Kan complexes in the terminology of \cite{kapranov1991infty}) in lemma $3.4$ (which would be sufficient for the rest of the argument) and to add in the definition of fibrant object a condition forcing cells to be ``weakly composable''. Note that this change in the definition of fibrant objects is probably necessary for the results of their section $2$ to be true without restricting the category of diagrams, as we do not see how they obtains a group structure on the $\pi_n$ without any assumption of this kind.

While this problem seems easily fixable, we believe that fixing it would probably only makes the ``real'' problems of this lemma $3.4$ appears: as we mentioned above, one cannot compose Johnson's diagrams in general, and that would provide other kind of counterexamples to this lemma, and secondly, degeneracies seems to allow to construct a presheaf such that an Eckmann-Hilton collapse happen in the free $\infty$-category it generates. In this case the identity cell going between some $u \#_0 v$ to a $v \#_0 u$ which are equal because of an Eckmann-Hilton collapse would not be representable by a single Johnson diagram with correct boundary. We initially wanted to give explicit examples of this two phenomenons, but unfortunately the first type of counterexamples related to the complete absence of compositions of cells whether the diagram compose or not appears so often that it seems nearly impossible to actually construct any interesting other kind of examples without first choosing a solution to this problem.
}

\bibliography{Biblio}{}
\bibliographystyle{plain}

\end{document}